

\documentstyle{amsppt}
\pageheight{20cm}
\pagewidth{16cm}

\NoBlackBoxes

\bf
\noindent THE THETA DIVISOR OF $SU_C(2,2d)^s$ IS VERY AMPLE IF $C$ IS
NOT HYPERELLIPTIC \footnote
{ The authors were partially supported by the European Science Project
"Geometry of Algebraic
Varieties" Contract no. SCI-0398-C(A).}
\rm  $$\text {by} $$    $$ \text {Sonia Brivio } \quad
\text {and} \quad \text {Alessandro Verra }  $$

\bigskip

\proclaim
{0.Introduction }
\endproclaim

\noindent
Let $C$ be an irreducible smooth complex curve
of genus $g \geq 2$: in this paper we are dealing with the moduli
space $$ X = SU_C(2,2d) \tag 0.0 $$
of semistable rank~2  vector bundles over $C$ having fixed  determinant
of even degree. As it is well
known  $$ X-Sing X = SU_C(2,2d)^s $$ is the moduli space for \it stable
\rm rank two vector bundles
and $Sing X$ is naturally isomorphic to the Kummer variety of the Jacobian of
$C$.  Let  $$ \Cal L \in
Pic X \tag 0.1 $$ be the generalized theta divisor $$ \theta :X \longrightarrow
\bold P H^0(\Cal L)^*
\tag 0.2 $$ the map associated to $\Cal L$, the main theorem we want to show is
the following
\proclaim {(0.3)THEOREM} Assume $C$ is not hyperelliptic, then \par \noindent
(1) $\theta$ is
injective \par \noindent (2) $d\theta_x$ is injective if $x \in X-Sing X$.
\endproclaim  \noindent
The proof of the injectivity of $d\theta_x$ at points $x \in Sing X$ seems a
more delicate
technical problem which is, at the moment, beyond the capability and perhaps
the patience of
the authors. Previous results on the embedding properties of $\theta$ were
obtained  by
Beauville [B1]:  \proclaim {(0.4)} \par \noindent (i) $\theta $
is a morphism of degree $\leq 2$ onto its image \par \noindent (ii) $ \theta $
has degree two if and
only if $C$ is hyperelliptic and $g \geq 3$ \endproclaim \noindent and by
Ramananan-Narasimhan for
$g \leq 3$ ([NR1], [NR2]). More recently  Laszlo showed in [L1] that \proclaim
{(0.5)}\par
\noindent $\theta $ is an embedding for a general non hyperelliptic curve $C$
\endproclaim \noindent
We recall that $Pic(X) \cong \bold Z$ and that by definition the generalized
theta divisor is the
ample generator of it; hence, in particular, $\theta $ is a finite morphism
over its image. The
fundamental geometric interpretation of $\theta$ is the following (cf. [B1]):
\par \noindent let $J =
Pic^0(C)$, $\Theta \subset J$ a symmetric theta divisor then, choosing as a
fixed determinant for the
points of $X$ the cotangent bundle $\omega_C$, there exists a canonical
isomorphism
 $$
 H^0(\Cal L)^* \cong H^0(\Cal O_J(2\Theta)) \tag 0.6
 $$
 Under the previous identification the map associated to $\Cal L$ becomes
 $$
 \theta : X \longrightarrow \mid 2 \Theta \mid \tag 0.7
 $$
\par \noindent on the other hand each semistable vector bundle $\xi $ having
determinant
$\omega_C$  defines a divisor
$$
\Theta_{\xi} \in \mid 2 \Theta \mid
$$
\par \noindent which is set theoretically so defined
$$
\Theta_{\xi} = \lbrace e \in J | h^0(\xi(e)) \geq 1 \rbrace \tag 0.8
$$
\par \noindent If $ x \in X$ is the moduli point of $\xi$ it turns out that
$$
\theta (x) = \Theta_{\xi} \tag 0.9
$$ \bigskip
\noindent Let us introduce very briefly the method we used to study $\theta$:
we have constructed a
\it projective \rm family $T$ of rational maps
$$
\theta_t: X \longrightarrow \bold P^n
$$
\par \noindent $n = \binom {g+1}2$, such that $\theta_t = \lambda_t \cdot
\theta$
\noindent where
$$
\lambda_t : \mid 2\Theta \mid \longrightarrow \bold P^n
$$
\noindent is a suitable linear projection depending on $t$. Globalizing this
construction, we have
a rational map
$$
F: X \times T \longrightarrow \bold P(\Cal Q)
$$
\noindent where $\Cal Q$ is a vector bundle on $T$ of rank $\binom {g+1}2 +1$
and $F/X\times t =
\theta_t$. It turns out that the maps $\theta_t$ have degree two onto their
image as well as the map
$F$. In particular the latter induces a birational involution
$$
j: X \times T \longrightarrow X \times T
$$
\noindent Of course $\theta$ is an embedding if any two  distinct points
(tangent vectors) of $X$
are separated by at least one $\theta_t$ ($d\theta_t$). Since $T$ is projective
we introduce in
section 6 a weaker version of rigidity lemma which is "ad hoc" for the rational
maps $F$ and
$j$. Applying this rigidity argument we finally deduce in section 7 that,
whenever $\theta$ does not
separates two given points of $X$ or two tangent vectors  of $X-Sing X$, then
the degree of $\theta$
must be two. Hence $C$ is hyperelliptic by Beauville's theorem.\par \noindent
Most of the paper is
necessarily devoted to a detailed study of the maps $\theta_t$, $F$ and $j$
(sections 2,3,4). This
seems interesting in itself in view of the geometry of quadrics behind
$\theta_t$. It is simple to
define $\theta_t$: the parameter space $T$ is actually $Pic^2(C)$ so that $t$
is a degree 2 line
bundle, consider the surface $$S_t = \lbrace t(-x-y), \quad x,y \in C \rbrace
\subset Pic^0(C)$$ then
$\bold P^n$ is just the linear system $\mid \Cal O_{S_t}(2\Theta) \mid$ and by
definition  $$
\lambda_t: \mid 2\Theta \mid \longrightarrow \bold P^n $$ is the restriction
map, so that $\theta_t =
\lambda_t \cdot \theta$. We will construct in section 4 a canonical projective
isomorphism $$
\mid \Cal O_{S_t}(2\Theta) \mid \cong \mid \Cal I_{C_t}(2) \mid
$$
where
$$
C_t \subset \bold P^{g+2}
$$
is the curve $C$ embedded by $\omega_C(2t)$ and $\Cal I_{C_t}$ its ideal sheaf.
This relates $X$
to the geometry of the quadrics containing $C_t$: the expected dimension for
the variety $Y$ of
quadrics of rank $\leq 6$  containing $C_t$ is $3g-3 = dim X$; let $Z_t =
\theta_t(X)$, we will show
that $Z_t$ is an irreducible component of $Y$. If $g \geq 3$ $\theta_t$ is the
"natural"
double covering of $Z_t$ which parametrizes the two rulings of maximal linear
subspaces in a rank 6
quadric, if $g=2$ $\theta_t = \theta$. The main technical problems to apply the
previous method are considered in section 5 where the birational involution $j$
is suitably extended
and  it is shown that $codim I(x) \geq 2$ for some special "Brill-Noether"
locus $I(x)
\subset T$ depending on $x \in X$: this follows essentially from Martens'
theorem [ACGH] and
Lange-Narasimhan results on Maruyama conjecture [LN]. \bigskip \noindent  This
paper originates in
part from the first author's doctoral thesis (Universit\'a di Torino, 1994) and
from our unpublished
manuscript [BV]. We wish to thank in particular A. Beauville, A.Collino,
I.Dolgachev, Y. Laszlo for
the interest and the helpful comments on this work. \bigskip \noindent Let us
fix here some of
the usual notations: \par \noindent $\xi$ $=$ semistable
rank two vector bundle on $C$ with determinant $\omega_C$, $[\xi]$ $=$ moduli
point of $\xi$ \par
\noindent $J$ $=$ Jacobian of $C$, $\Theta$ $=$ a symmetric theta divisor in
$J$, $C(n)$ = n-th
symmetric product of $C$\par \noindent $T$ $=$ $Pic^2(C)$, $t$ $=$  an element
of $T$ \par \noindent
$\bold P^{g+2}_t = \bold PH^0(\omega_C(2t))^*$, $C_t$ $=$ the curve $C$
embedded in $\bold P^{g+2}_t$
by $\omega_C(2t)$. \par \noindent $E^*$ $=$ dual of the vector space (bundle)
$E$. \par \noindent
\proclaim {1. Preliminaries: rank 6
quadrics and rank 2 vector bundles} \endproclaim \noindent For our and reader's
convenience we
recollect in this section  some standard properties  which will be used many
times. We fix  a smooth,
irreducible projective curve  $$ C \subseteq \bold P^n $$
which is linearly normal and not degenerate. The space of global sections of
$\Cal
O_C(1)$ will be
$$
H = H^0(\Cal O_C(1)) \tag 1.1
$$
so that $\bold P^n = \bold PH^*$. Then we consider a pair
$$
(\Cal E, V) \tag 1.2
$$
such that: \proclaim {(1.3)}\par \noindent
(i) $\Cal E$ is a rank 2 vector bundle over $C$  \par
\noindent
(ii) $V$ is a 4-dimensional vector space in $ H^0(\Cal E) $ \par \noindent
(iii) $det (\Cal E) = \Cal O_C(1)$ \par \noindent
(iv) $h^0(\Cal E^*) = 0$
\endproclaim \noindent  from the previous data
we construct:  \par \noindent \proclaim {(1.4)} \par \noindent (i) the
evaluation map $ e_V: V \otimes \Cal O_C \rightarrow \Cal E $ \par \noindent
(ii) the determinant map $ d_V: \wedge ^2 V \rightarrow H $. \par \noindent
(iii) the Grassmannian $G^*_V$ of 2-dimensional subspaces of $V^*$
\endproclaim   \noindent
By definition $G^*_V$ is the projectivized set of the irreducible vectors in
$\wedge^2 V^*$ hence
$$G^*_V \subset \bold P^5 = \bold P(\wedge^2 V^*) $$ as a smooth quadric
hypersurface. The dual of
of the evaluation map ${e_V}^*:{\Cal E }^* \rightarrow \Cal O_C \otimes V^*$
defines in $V^*$ the
family of subspaces $$ \lbrace V^*_x = Im e_{V,x}^*, x \in C \rbrace$$
\noindent By definition $e_V$
is \it generically surjective \rm iff $dim V^*_x = 2$ for a general $x$,
assuming this one can
construct a rational map  $$ g_V: C \longrightarrow G_V^* $$ \noindent which
associates to $x$ the
point $$g_V(x) = \wedge ^2 V^*_x \in G^*_V$$ \proclaim {(1.5) DEFINITION} $g_V$
is the Gauss map of
the pair $(\Cal E,V)$ \endproclaim \proclaim {(1.6)PROPOSITION} \par \noindent
(1) $e_V$ is
generically surjective if and only if the determinant map  $d_V$ is not zero
\par \noindent (2)
assume $d_V$ is not zero,  let $ \delta_V : \bold P^n \longrightarrow \bold P
\wedge^2V^* $ be the
projectivization of the dual map $d_V^*$, then $$ \delta_V /C  =  g_V  $$
\noindent (3) $\delta_V$
is defined at $x \in C$ iff $e_{V,x}$ is surjective, in particular $Supp Coker
(e_V) \cong C\cdot
\bold P Ker({d_V}^*) $   \endproclaim \demo {Proof}ÊWell known. \enddemo
\noindent Let  $p $ be an equation for the quadric $G^*_V$  and let $$ q(\Cal
E,V)  $$ be  the
pull back of $p$ by the dual map $d^*_V : H^* \rightarrow \wedge^2 V^*$. If it
is not identically zero
$q(\Cal E,V)$ defines a quadric in $\bold P^n$. This will be denoted by $$
Q(\Cal E,V)
$$
\par \noindent and it is the most important object of this paper:
\proclaim{(1.7) DEFINITION} $q(\Cal E,V)$ is a quadratic form defined by $(\Cal
E,V)$. $Q(\Cal E,V)$
is the quadric of $(\Cal E,V)$. \endproclaim \noindent  Let $ q = q(\Cal E,V)$,
$ Q = Q(\Cal
E,V) $ obviously $Q =\delta_V^{-1}(G^*_V)$  so that
$$
 \text { rank $q \leq  6$ and $q$ vanishes on $C$} \tag 1.8
$$
\par \noindent ÊLet $K = Ker ({d_V}^*)$, $I = Im d_V^*$ then $Q$ can be
considered as a cone of vertex
$\bold PK$ over the quadric $\bold PI \cap G^*_V$. In particular $\bold PK
\subseteq Sing Q$. We want
to point out that \bigskip \noindent (1.9) (i) $ \bold PK = Sing Q $ if and
only if $\bold PI$ is
transversal to $Q$ \par \noindent (ii) if $q$ has rank $\geq 5$ then $ Sing Q =
\bold PK$ \bigskip
\noindent \par \noindent For any subline bundle $L \subset \Cal E$ we have in
$H^0(\Cal E)$ the
natural intersection  $$ V_L = V \cap H^0(L)  $$ \proclaim {(1.11) PROPOSITION}
Let $q = q(\Cal
E,V)$, $r$ the rank of $q$. Then: \par \noindent (1) $r \leq 4
\Longleftrightarrow $ $\Cal E$
contains a a subline bundle $L$ such that $dim V_L \geq 2$ \par \noindent (2)
$r = 0
\Longleftrightarrow \Cal E$ contains a subline bundle $L$  such that $dim V_L
\geq 3$ \par \noindent
(3) $e_V$  not generically surjective $\Longleftrightarrow \Cal E$ contains a
subline bundle $L$
such that $dim V_L =
 4$. \endproclaim  \demo {Proof} By duality the elements of $\wedge ^2 V$ are
linear forms on
$\wedge^2 V^*$. Let $G_V \subset \bold P (\wedge^2 V)$ be the projectivized set
of the irreducible
vectors: $G_V$ is the dual of the  quadric  $G_V^*$, hence a non zero
irreducible vector $s_1
\wedge s_2 \in \wedge^2 V$ is a linear form defining a singular hyperplane
section of $G_V^*$. By
duality again $Im(d_V^*)$ is the zero locus of the linear forms in $Ker d_V$.
{}From this it follows
that: \par \noindent (a) $r \leq 4$ $\Longleftrightarrow$ $Ker d_V$ contains a
non zero irreducible
vector $s_1 \wedge s_2$ \par \noindent (b) $r = 0$ $\Longleftrightarrow$  $Ker
d_V$ contains a
3-dimensional vector space $F$ of irreducible vectors. \par \noindent
On the other hand it is a standard property of two independent sections
$s_1,s_2 \in V$ that \par
\noindent (c) $d_V(s_1 \wedge s_2) =0$ $\Longleftrightarrow$ $s_1,s_2$ define a
subline bundle $L
\subset \Cal E$ such that $s_1,s_2 \in V_L$. \par \noindent
The statement follows easily from (a), (b), (c): we omit the details for
brevity.\enddemo
\bigskip \noindent  We  will say that two pairs $(\Cal E_1,V_1)$,
$(\Cal E_2,V_2)$  are \it isomorphic \rm if there exists an isomorphism $\sigma
: \Cal E_1
\rightarrow \Cal E_2$ such that $\sigma^* (V_2) = V_1$.Ê We want to construct
all the isomorphism classes of pairs $(\Cal E,V)$ as in 1.3 which define the
same quadratic form $q$:
we are  interested to  do this only when $q$ has rank $\geq 5$ and $Sing Q \cap
C = \emptyset$, ($Q$
$=$ $\lbrace q = 0 \rbrace $). \bigskip \noindent \proclaim {(1.12)
DEFINITION} $V$ \it
generates \rm $\Cal E$ iff $Coker(e_V)$ = $0$. \endproclaim \noindent Assume
$V$ generates $\Cal E$
then in the exact sequence $$ \CD 0 @>>> {\Cal E^*} @>{e_V^*}>> {V^* \otimes
\Cal O_C} @>>> {Coker(e_V^*)} @>>> 0 \\ \endCD \tag 1.13 $$ $Coker(e_V^*)$ is
locally free and $V^*$
is canonically identified to a 4-dimensional space of global sections of it.
This defines another
pair $$ (\overline {\Cal E}, \overline V) = (Coker(e_V^*),V^*)  $$ \noindent
satisfying the
assumptions in (1.3). Repeating the same construction for $(\overline {\Cal E},
\overline V)$ one
obtains $(\Cal E,V)$ again.\proclaim {(1.14) DEFINITION} $(\overline {\Cal E},
\overline V)$ is the
\it dual pair \rm of $(\Cal E,V)$ \endproclaim \noindent The \it universal \rm
and the \it quotient
 bundle \rm of $G^*_V$ will be denoted respectively by $$\Cal U_V \quad \text
{,} \quad \overline
{\Cal U}_V \tag 1.15$$ \noindent from the definition of the Gauss map $g_V$ it
follows   $$  g_V^*
\Cal U_V \cong Im(e_V^*) \cong Im(e_V)^* \quad \text {,} \quad g_V^* \overline
{\Cal U}_V \cong
Ker(e_V^*) \cong Ker(e_V)^* \tag 1.16$$ and the exact sequence $$ 0
\longrightarrow g^*_V \Cal U^*_V
\longrightarrow \Cal E \longrightarrow Coker (e_V) \longrightarrow 0
 $$ \noindent In particular $$ \Cal E \cong g_V^* \Cal U_V^*\quad \text{,}
\quad
\overline {\Cal E} \cong g_V^* \overline {\Cal U}_V \tag 1.17 $$ if and only if
$V$ generates $\Cal
E$. In this case the previous sequence (1.13) is just the pull-back by $g_V$ of
the standard
universal/quotient bundle sequence $$
0 \longrightarrow \Cal U_V \longrightarrow V^*\otimes \Cal O_{G^*_V}
\longrightarrow
\overline {\Cal U}_V \longrightarrow 0$$ \proclaim {(1.18) LEMMA} Let $(\Cal
E_1,V_1)$, $(\Cal
E_2,V_2)$ be two pairs such that \par \noindent (i) $q = q(\Cal E_1,V_1) =
q(\Cal E_2,V_2)$
\par \noindent (ii) $V_i$ generates $\Cal E_i$, $i = 1,2$ \par \noindent (iii)
$Ker(d_{V_1}^*) =
Ker(d_{V_2}^*)$ \par \noindent Then either $(\Cal E_1,V_1)$, $(\Cal E_2,V_2)$
are isomorphic or they
are dual. Moreover: \par \noindent
let $h^0(\Cal E_1) = 4$, if the two pairs are dual and $\Cal E_1 \cong \Cal
E_2$ then $q$ has rank
$\leq 5$.\endproclaim \demo {Proof} For $i=1,2$ consider the Gauss maps $
g_{V_i}: C \longrightarrow
G_{V_i}^*  $: since $V_i$ generates $\Cal E_i$ it follows  $\Cal E_i^* \cong
g_{V_i}^* \Cal U_{V_i}$.
Since $Ker(d_{V_1}^*)$ = $Ker(d_{V_2}^*)$  there exists a unique isomorphism
$ \phi : Im(d_{V_1}^*)
\rightarrow Im(d_{V_2}^*)$ such that $d_{V_2}^* = d_{V_1}^* \cdot \phi $. Let
$v \in \wedge^2 V^*_1$
be an irreducible vector, $v = d_{V_1}^*(h)$ where $h \in H^*$ then $q(\Cal
E_1,V_1)$ is zero on $h$.
By (i) the same holds for $q(\Cal E_2,V_2)$  so that $d_{V_2}^*(h)$ $=$ $\phi
(v) $ is an irreducible
vector too. Since it preserves irreducible vectors  $\phi$ is the restriction
to $Im d_{V_1}^*$ of an
isomorphism $\Phi : \wedge ^2 V_1^* \rightarrow \wedge^2 V_2^*$ with the same
property. Therefore
$\Phi$ induces  a biregular isomorphism $$ \Phi/G_{V_1}^* = f:
G_{V_1}^*Ê\longrightarrow G_{V_2}^*
\quad \text { such that}Ê\quad g_{V_2}=f \cdot g_{V_1}$$ \noindent  As it is
well known either $\Cal
U_{V_1} \cong f^* \Cal U_{V_2}$ or $f^*\Cal U_{V_2} \cong \overline {\Cal
U}_{V_1}^*$.
Assuming $f^* \Cal U_{V_2}$ $ \cong $ $\overline {\Cal U}_{V_1}^*$ and applying
$f^*$ to the
universal/quotient bundle sequence   $$ 0 \rightarrow \Cal U_{V_2} \rightarrow
V_2^* \otimes
\Cal O_{G_{V_2}} \rightarrow \overline{\Cal U}_{V_2} \longrightarrow 0
 $$ \noindent one has  $$ 0
\rightarrow \overline {\Cal U}_{V_1}^* \rightarrow V_1 \otimes
\Cal O_{G_{V_1}} \rightarrow \Cal U_{V_1}^* \rightarrow 0
$$
 \noindent pulling this sequence back by  $g_{V_1}$ and using
$g_{V_1}^*$ $\cdot $ $f^*$ = $g_{V_2}^*$ we finally obtain  $$ 0Ê\rightarrow
\Cal E_2^* \rightarrow V_1 \otimes \Cal O_C \rightarrow \Cal E_1
\rightarrow 0$$ \noindent  This precisely means that the two pairs are
dual. If $f^* \Cal U_{V_2} \cong \Cal U_{V_1}$ the same argument implies
that the two pairs are isomorphic. \noindent To prove the latter statement
observe that, by
assumption, we can set $\Cal E_1 = \Cal E_2 = \Cal E$ and $V_1 = V_2 = V$.
Then, tensoring  the
previous sequence by $\Cal E$ and passing to the long exact sequence, we obtain
$$
0 \rightarrow H^0(\Cal E \otimes \Cal E^*) \rightarrow V \otimes V \rightarrow
H^0(\Cal E \otimes
\Cal E)
$$
\noindent where the latter arrow is the natural map $\mu$. Let $\mu^-$ be the
restriction of $\mu$ to
$\wedge^2V \subset V \otimes V$, then $Im \mu^-$ is contained in $H^0(det \Cal
E) \subset H^0(\Cal E
\otimes \Cal E)$ and moreover $\mu^-$ is just the determinant map $d_V$. Since
$h^0(\Cal E \otimes
\Cal E^*) \geq 1$ $\mu$ is not injective: we leave as an exercise to check that
then $\mu^- = d_V$
must be not injective. Since the rank $r$ of $q$ is the dimension of $Im d_V$
it follows $r \leq
5$.\enddemo  \proclaim {(1.19)PROPOSITION} Let $Q$ be a quadric of rank $r \leq
6$ which contains
$C$, then there exists a pair $(\Cal E,V)$ such that $Q = Q(\Cal E,V)$. If $r$
= $5,6$ and  $Sing Q
\cap C = \emptyset$  every pair defining $Q$ is isomorphic to $(\Cal E,V)$ or
to its dual $(\overline
{\Cal E}, \overline V)$. If the rank is 5 these pairs are isomorphic too.
\endproclaim \demo {Proof}
Consider the linear projection $\delta : Q  \rightarrow \bold P^{r-1}$ of
center $Sing Q$ and the
quadric $\delta(Q)$: $\delta(Q)$ can be considered as a linear section of $G$,
where $G \subset \bold
P^5$ is the Pluecker embedding of a Grassmannian $G(2,4)$. Let $\Cal U$ be the
universal bundle on
$G$, $g = \delta / C$,  since $Sing Q \cap C = \emptyset$ the pull-back $g^*
\Cal U^*$ has
determinant $\Cal O_C(1)$. Putting $$V = g^*H^0(\Cal U^*) \quad , \quad \Cal E
= g^*\Cal U^*$$  we
obtain a pair $(\Cal E,V)$ as in (1.3). It is easy to see that this pair
defines $Q$, assume $ r \geq
5$ then $\bold PKer (d_V^*)$ = $ Sing Q$ (see 1.9) and $V$ generates $\Cal E$
because $Sing Q \cap C
= \emptyset$ (see 1.6). Note that these properties hold for every pair defining
$Q$: let $(\Cal
E_1,V_1)$, $(\Cal E_2,V_2)$ be two of them, then, by the previous lemma 1.18,
either $(\Cal
E_1,V_1)$, $(\Cal E_2,V_2)$ are isomorphic or they are dual. Finally let $r=5$
then $\delta(Q)$ is a
smooth hyperplane section of $G$ and hence the universal and the quotient
bundle restrict to the same
bundle on $\delta(Q)$ (see e.g.[O]). This easily implies that $(\Cal E,V)$ is
isomorphic to its dual
and completes the proof.  \enddemo \noindent Now we fix a rank two vector
bundle $\Cal E$ with $det
\Cal E = \Cal O_C(1)$ and the space of its global sections $$ W = H^0(\Cal E)
\tag 1.20 $$ we assume
$dim W \geq 4$. Let $\Cal I_C$ be the ideal sheaf of $C$,  the construction of
$q(\Cal E,V)$ defines
a linear map $$ h_{\Cal E}: \wedge^4 W \longrightarrow H^0(\Cal I_C(2)) \tag
1.21$$  which is so
defined on irreducible vectors:  let $v \in \wedge^4 W$, $V$  the
corresponding   4-dimensional
subspace of $W$ then $$ h_{\Cal E}(v) = q(\Cal E,V) $$ We are more interested
to the following
specialization of the of the previous map: Let us fix a 3-dimensional vector
space $$ F \subset W
\tag 1.22 $$ and the image $$ H_F \subseteq H \tag 1.23 $$ of $\wedge^2 F$
under the determinant map
$$ d: \wedge^2W \rightarrow H = H^0(\Cal O_{\bold P^n}(1))$$ we restrict
$h_{\Cal E}$ to $$ W_F =
(\wedge^3F) \wedge W \subseteq \wedge^4 W $$ Every non zero $v \in W_F$
represents a 4-dimensional
vector space $V$ which contains $F$, in particular $q(\Cal E,V)$
 vanishes on the linear space  $ \Lambda_F \subset \bold P^n$ which is the
zero locus of the elements of $H_F$ and has codimension $\leq 3$ in $\bold
P^n$.
Let $$ \Cal I_F = \quad \text {  ideal sheaf of $C \cup \Lambda_F$} $$ then
$h_{\Cal E}$ restricts to a linear map $$ h_F: W_F \longrightarrow H^0(\Cal
I_F(2)) \tag 1.24 $$ \proclaim {(1.25)PROPOSITION} Let $e_F : F \otimes \Cal
O_C
\rightarrow \Cal E$ be the evaluation map, $d_F : \wedge^2 F \rightarrow H$ the
determinant map: \par \noindent (1) if $d_F$ is injective then $h_F$ is
injective \par \noindent (2) if $Coker (e_F) = 0$  then $h_F$ is an
isomorphism.
\endproclaim \demo {Proof} (1) Let $F \subset V$, $v = \wedge^4V$, $q = q(\Cal
E,V)$. By proposition 1.11 $q = h_F(v)$ is zero iff $\Cal E$ contains a subline
bundle $L$ such that $dim V_L \geq 3$. Assume such an $L$ exists and consider
$F_L = V_L \cap F$: then $dim F_L \geq 2$; moreover, by remark (c) in the proof
of proposition 1.11, the determinant map is zero on  $\wedge^2 V_L$ hence on
$\wedge^2 F_L$: against the injectivity of $d_F$.  \par \noindent (2) Consider
the dual map  $ e_F^*: \Cal E^* \rightarrow \Cal O_Y \otimes F^*$ and in
$F^*$ the family of subspaces  $ F^*_x = Im(e_{F,x}^*)$, $x \in C
$. Since $Coker (e_F) = 0$ $dim F^*_x = 2$ for all $x \in C $. Then the  Gauss
map $$ g_F: C \longrightarrow \bold P^2 = \bold P \wedge^2 F^*  $$ sending $x$
to
$\wedge^2 Im(e_{V,x}^*)$ is a morphism and hence $$ \Cal E^* \cong
g_F^* \Cal \Omega_{\bold P^2}(1) $$ On the other hand $g_F$ is  defined by the
subspace $H_F$ of $H = H^0(\Cal O_C(1))$ which means that $g_F$ is the
restriction to $C$ of the linear projection of center $\Lambda_F$. To
show the surjectivity of $h_F$ consider $q' \in H^0(\Cal I_F(2))$: since it
vanishes on $\Lambda_F$  $q'$ has rank $\leq 6$. Hence $q'$ = $q(\Cal E',V')$
for a given pair $(\Cal E',V')$ (prop. 1.19). Up to replacing this pair by its
dual we can assume that $V'$ contains a 3-dimensional subspace $F'$ such that
the
image of $\wedge^2 F'$ by the determinant map is again $H_F$. For this reason,
when we construct as above the Gauss map $$ g_{F'}: C \longrightarrow \bold P^2
$$ from $F'$, we obtain $g_F$ = $g_{F'}$ = linear projection of center
$\Lambda_F$. Therefore   $\Cal E' \cong \Cal E$ and  hence $q' \in Im (h_F)$.
By
(1), to show the injectivity of $h_F$ it suffices to show that $d_F$ is
injective. Assume $d_F(s_1 \wedge s_2) = 0$ for a given $s_1 \wedge s_2 \in
\wedge ^2 F$ then, for any $V$ containing $F$, $q(\Cal E,V)$ has rank $\leq 4$
(prop. 1.11). This immediately implies $\Lambda_F \cap C \neq \emptyset$ and
$Coker (e_F) \neq 0$.\enddemo
\proclaim
{(1.26)PROPOSITION}ÊAssume $\Lambda_F \cap C = \emptyset$.
Then $h_F$ is an isomorphism and in particular $$ h^0(\Cal E) = 3 + h^0(\Cal
I_F(2)) $$ \endproclaim \demo {Proof} Observe that $\Lambda_F \cap C =
\emptyset$ implies that $e_F$ is everywhere surjective. Hence, by the previous
proposition, $h_F$ is an isomorphism. In particular $h^0(\Cal I_F((2))  = dim
W_F = dim W - 3$ \enddemo \bigskip

\proclaim
{\bf 2. The fundamental involution and the fundamental map on $SU_C(2,2d)
\times Pic^2(C)$ \rm} \endproclaim \bigskip \noindent Let $$ U_C(2,2d) \tag 2.1
$$ be
the moduli space of semistable rank 2 vector bundles $\Cal E$ on $C$ of fixed
degree $2d$. Since $U_C(2,2d)$ does not depend on the choice of $d$ there is no
restriction to assume $$ 2d = 2g+2  $$\noindent then the construction of the
dual pair we
have given in the previous section defines an involution  $$  i: U_C(2,2d)
\longrightarrow U_C(2,2d) \tag 2.2 $$ Roughly speaking $i$ is so defined: let
$u
\in U_C(2,2d)$ be the moduli point of $\Cal E$, since $deg(\Cal E) = 2g+2$  the
expected dimension of $V$ = $H^0(\Cal E)$ is 4. Assuming this and that $V$
generates $\Cal E$ we can construct  the dual $(\overline {\Cal E}, \overline
V)$ of the pair $(\Cal E,V)$. If  semistable $\overline {\Cal E}$ defines
another point $\overline u \in U_C(2,2d)$,  we  define in this case $$i(u) =
\overline u$$\noindent Let  $$ X = SU_C(2,2d) \quad T = Pic^2(C) \tag 2.3
$$\noindent for  reasons
which will be evident later we are more interested to  the analogous  $$ j: X
\times T \longrightarrow X \times T \tag 2.4 $$\noindent of the involution $i$:
let
$(z,t) \in X \times T$, then $z$ is the moduli point of $\xi$ where $det(\xi) =
\omega_C$ and $t$ is a line bundle of degree 2. Putting $\Cal E = \xi(t)$, $V$
=
$H^0(\xi(t))$ we obtain a pair $(\Cal E,V)$ as above. Passing to the dual pair
we obtain a rank 2 vector bundle $\overline {\Cal E}$ such that $det (\overline
{\Cal E}(-t))$ = $\omega_C$. If semistable $\overline {\Cal E}(-t)$   defines a
point $\overline z \in X$, then we define as above $$ j(z,t) = (\overline z,t)
$$ The relation between $i$ and $j$ is very clear:  let $$ m: X \times T
\longrightarrow U_C(2,2d) \tag 2.5 $$\noindent be the  tensor product map
sending  the
moduli  of the pair $(\xi,t)$ to the moduli of the vector bundle $\Cal E =
\xi(t)$, from the definitions of $i$ and $j$ we have the commutative diagram
$$
\CD {SU_C(2,2d) \times T} @>>j> {SU_C(2,2d) \times T} \\ @VVmV
@VVmV \\
{U_C(2,2d)}Ê@>>i> {U_C(2,2d)} \\ \endCD \tag 2.6 $$ It is not difficult to see
that
$m$ is an \'etale covering of degree $2^{2g}$, in particular  $j$ can be viewed
as the lifting of $i$ to $X \times T$. To begin we construct  a suitable open
subset $$U_j \subset X
\times T \tag 2.7 $$  where $j$ will be properly defined and biregular. $U_j$
will be useful
later, to define it we consider the exact sequence $$ 0 \rightarrow \Omega
\rightarrow H^0(\omega_C) \otimes \Cal O_C \rightarrow \omega_C \rightarrow 0
\tag 2.8 $$ $\Omega $ is a well known semistable vector bundle of rank $g-1$
over $C$. With some abuse we will use sometimes the same notation for a
semistable bundle and for its moduli point: \proclaim {(2.9) DEFINITION} $U_j =
\lbrace (\xi,t) \in X \times T / h^0(\Omega \otimes \xi(t)) = 0 \rbrace $
\endproclaim  \par \noindent We want to point out that , tensoring 2.8  with
$\xi(t)$ and passing to the long exact sequence, it follows $$ h^0(\Omega
\otimes \xi(t))=0 \Leftrightarrow \text {the multiplication map $ \nu
:H^0(\omega_C) \otimes H^0(\xi(t)) \rightarrow H^0(\omega_C \otimes \xi(t))$ is
injective } \tag 2.10 $$ furthermore, since $\chi(\Omega \otimes \xi(t)) = 0$,
$$ h^0(\Omega \otimes \xi(t)) = 0 \Leftrightarrow \nu : H^0(\omega_C) \otimes
H^0(\xi(t)) \rightarrow H^0(\omega \otimes \xi(t)) \text {is an isomorphism}
 $$ this will be used various times.  Of course $U_j$ is open,
\proclaim {(2.11)PROPOSITION}Let $(z,t) \in Sing X \times T$, then \par
\noindent
(i) $U_j \cap Sing X \times  t  \neq \emptyset $ so that $U_j \cap X \times t$
is not empty \par
\noindent
(ii) $ U_j \cap \lbrace z \rbrace \times T \neq \emptyset$. \endproclaim \demo
{Proof}
Let $z \in Sing X$,  it is well known that $z$ is the moduli point of all
semistable extensions  $$ 0
\rightarrow \alpha \rightarrow \xi \rightarrow \omega_C \otimes \alpha^{-1}
\rightarrow 0 \quad \text {or} \quad  0 \rightarrow \omega_C \otimes
\alpha^{-1}
\rightarrow \xi \rightarrow \alpha \rightarrow 0$$ for a given $\alpha \in
Pic^{g-1}(C)$. Tensoring the previous sequences by $\Omega(t)$ it follows that
$h^0(\Omega \otimes \xi(t))=0$ iff $h^0(\Omega \otimes \alpha(t)) $ =
$h^0(\Omega \otimes \omega \otimes \alpha^{-1})$ = $0$.  Tensoring the sequence
2.8 by $\alpha(t)$  the latter condition is equivalent to the injectivity of
the
multiplication  $ H^0(\omega_C) \otimes H^0(\alpha(t)) \longrightarrow
H^0(\omega_C \otimes \alpha(t)) $ and of the analogous multiplication for
$\omega_C \otimes \alpha^{-1}$. Fixing $t$ the injectivity of these maps for a
general $\alpha$ is well known and it follows from the base-point-free pencil
trick [ACGH]. Fixing $\alpha$ (hence $z$) the same argument gives the
injectivity
for a general $t$. Hence $ U_j$ intersects both $z \times T$ and $Sing (X)
\times t$. \proclaim {(2.12) PROPOSITION } $(\xi,t)$ defines a point of $U_j$
if
and only if the following conditions are satisfied: \par \noindent (1)$
h^0(\xi(t))=4$\par \noindent (2) $\xi(t)$ is globally generated \par \noindent
(3) $h^0(\overline {\xi}(t))=4$, where $\overline {\xi}$ is defined by the
fundamental exact sequence $$ 0 \longrightarrow (\xi(t))^* \longrightarrow
H^0(\xi(t))^*\otimes \Cal O_C \longrightarrow \overline {\xi}(t)
\longrightarrow
0 \tag 2.13 $$ \endproclaim \demo {Proof} If (1),(2),(3) hold  the sheaf
$\overline {\xi}(t)$ defined by the previous sequence is a rank 2 vector bundle
such that   $h^1(\overline {\xi}(t)) = 0$. The semistability of $\xi(t)^*$
implies $h^0(\xi(t)^*) = 0$ so that, taking the  long exact sequence of 2.13,
we
obtain the isomorphisms $$ H^0(\xi(t))^* \cong H^0(\overline {\xi}(t)) \quad
\text {and} \quad H^1(\xi(t)^*) \cong H^0(\xi(t))^* \otimes H^1(\Cal O_C))  $$
By
Serre duality   the latter one is the dual of  the multiplication map $$\nu:
H^0(\omega_C) \otimes H^0(\xi(t)) \rightarrow H^0(\omega \otimes \xi(t))$$
since
$\nu$ is an isomorphism $h^0(\Omega \otimes \xi(t)) $ = $0$ and   $(\xi,t)$ $
\in $ $U_j$. Conversely let $\nu$ be an isomorphism: note that  $h^1(\omega_C
\otimes \xi(t))$ = $h^0(\xi(t)^*)$; then $h^0(\omega_C \otimes \xi(t)) = 4g$
and
$h^0(\xi(t)) \leq 4$. Hence, by Riemann-Roch, $h^0(\xi(t)) = 4$ and  (1) holds.
{}From the semistability of $\xi(t-x)^*$ we have again $ 0 = h^0(\xi(t-x))^*$ =
$h^1(\omega_C \otimes \xi(t-x))$, $\forall x \in C$. In particular it follows
$$h^0(\omega_C \otimes \xi(t-x)) = 4g-2 \quad \forall x \in C$$ so that
$\omega_C \otimes \xi(t)$  globally generated. Let $x \in C$,  consider in
$H^0(\omega_C) \otimes H^0(\xi(t))$ the vector space $ S_x = \nu ^{-1}(
H^0(\omega_C \otimes \xi(t-x)) $: this is just the sum $$ S_x =
H^0(\omega_C(-x)) \otimes H^0(\xi(t)) + H^0(\omega_C) \otimes H^0(\xi(t-x))$$
moreover the intersection of the latter two vector spaces is $H^0(\omega_C(-x)
\otimes H^0(\xi(t-x))$. If $\xi(t)$ is not globally generated at $x$ on sees by
computing dimensions that $dim S_x > 4g-2 = dim \nu(S_x)$: against the
injectivity of $\nu$. Hence (2) holds. Since (1) and (2) hold it is very easy
to
show (3) using  the  exact sequence  2.13.  \enddemo  \proclaim
{(2.14)PROPOSITION } Assume that $(\xi,t)$ defines a point $U_j$. Let
$\overline
{\xi}$ be constructed from $(\xi,t)$ using the exact sequence 2.13. Then
$\overline{\xi} $ is semistable and  $(\overline {\xi},t)$ defines a point of
$U_j$. \endproclaim \demo {Proof} Assume $\overline {\xi}$ not semistable the
there exists an exact sequence of vector bundles $$ 0 \longrightarrow L
\longrightarrow \overline {\xi}(t) \longrightarrow M \longrightarrow 0 $$ such
that $deg(L) \geq g+2$. By the previous proposition we have
$h^0(\overline{\xi}(t)) = 4$, therefore, passing to the long exact sequence, we
obtain $h^1(\overline {\xi}(t)) = h^1(M)=0$. Then only two  cases are
possible:   (1) $h^0(L)=3$ so that $deg(M) = g$, $h^0(M)=1$ and (2) $h^0(L) =
4$. In case (1) $\overline {\xi}(t)$ is not globally generated at points of
$Supp M$, in case (2) at all points of $C$. On the other hand $\overline
{\xi}(t)$ is defined by the exact sequence 2.13 hence it is globally generated.
Therefore $\overline {\xi}$ must be semistable. The rest of the proof is an
immediate consequence of prop. 2.12. \enddemo \proclaim {(2.15)PROPOSITION}
Assume that $(\xi,t)$ defines a point of $U_j$, then there is no subline bundle
$L$ of $\xi(t)$ having degree $\leq g$ and $h^0(L) \geq 2$. \endproclaim \demo
{Proof} As above we construct from $(\xi,t)$ the pair $(\overline {\xi},t)$. We
 know from the previous proposition that $\overline {\xi}$ is semistable.
Assume such a line bundle $L$ exists, there is no restriction to assume that
$L$
is saturated in $\xi(t)$. Let    $h^0(L)=2$, $L$ globally generated,   $deg(L)
\leq g$. Then we have  the exact diagram $$ \CD @.       0   @.
0                   @.   0  @. @. \\ @.       @VVV
@VVV                      @VVV    @. \\ 0 @>>> {L^*} @>>> {H^0(L)\otimes \Cal
O_C } @>>> L @>>> 0  \\ @. @VVV  @VVV  @VVV  @. \\ 0 @>>> {\overline
{\xi}(t)^*}
@>>> {H^0(\xi(t))Ê\otimes \Cal O_C } @>>> {\xi(t)}Ê@>>> 0 \\ @. @VVV @VVV @VVV
@. \\ 0 @>>> {M^*} @>>> {H^0(M) \otimes \Cal O_C}Ê@>>> M @>>> 0 \\ @. @VVV @VVV
@VVV @. \\ @. 0 @. 0 @. 0 @. @. \\Ê\endCD $$\noindent  where the middle
horizontal row is
just the dual of 2.13. In particular $L^*$ is a subline bundle of $\overline
{\xi}(t)^*$ and $deg(L^*) > \frac 12 deg(\overline {\xi}(t)^*)$: against the
semistability of $\overline {\xi}(t)$.   Assume $h^0(L) \geq 3$ or $h^0(L)=2$
and $L$ not globally generated then it follows $h^0(\xi(t-x)) \geq 3$ for at
least one $x \in C$. This implies $h^0(\xi(t)) > 4$  or $\xi(t)$ is not
globally
generated: a contradiction. \enddemo \noindent Using the exact sequence 2.13
we   define a bijective involution $$ j: U_j \longrightarrow U_j $$\noindent
such that $$
j(\xi,t) = (\overline {\xi},t) $$\noindent it is standard to check that  $j$ is
a
birational involution on $X \times T$ which is biregular on $U_j$. \proclaim
{(2.16)DEFINITION} Let $X = SU_C(2)$, $T = Pic^2(C)$. We will say that $$ j: X
\times T \longrightarrow X \times T $$\noindent is the fundamental involution
of $X
\times T$. \endproclaim
\noindent (2.17)REMARKS \par \noindent
(1) Let $K$ be the kernel of the evaluation map $e: H^0(\xi(t)) \otimes \Cal
O_C \rightarrow
\xi(t)$: it is standard to deduce that $j(\xi,t) = (K^*(-t),t)$ as soon as $K$
is a rank two vector
bundle which is semistable and of degree $2g+2$. It is not difficult to check
for which $t$ $K$ is
a rank two vector bundle and it is possible to suitably extend $j$ in case
$\xi(t)$ is not globally
generated (see section 5). The main problem is the stability of $K$, (cfr.
[Bu]). \par
\noindent
(2)Let $p_1: X \times T \rightarrow T$ be the first projection: knowing that $X
= p_1(U_j)$
would be very useful in this paper. This property means that, for each
semistable rank 2 vector
bundle $\xi$ with canonical determinant, the locus $$ D_{\xi} = \lbrace t \in T
/ h^0(\Omega \otimes
\xi(t)) > 0 \rbrace $$ is a divisor and not all $T$. In any case let $Y =
X-p_1(U_j)$: $Y$ is clearly
closed and contained in $X-Sing (X)$ by prop. 2.11. Since $Pic (X) \cong \bold
Z$ $Y$ cannot be a
divisor: $$  codim_X Y \geq 2 \tag 2.18 $$This remark will be important in the
proof of the main
theorem (see section 7 and thm. 6.13).
  \noindent We continue with the previous
notations and introduce in addition the following ones:  \par
\noindent (2.19) let $t \in T$ \par \noindent $ \bold P^{g+2}_t = \bold
PH^0(\omega_C(2t))^*   $ \par
\noindent $ C_t$ =  the curve $C$ embedded in $\bold P^{g+2}_t$ by
$\omega_C(2t)$    $$ \Cal I_t \text {= ideal sheaf of $C_t$}Ê \tag 2.20  $$ We
want to construct a vector bundle $$ \Cal Q \longrightarrow T $$ of fibre
$H^0(\Cal I_t(2))$ over $t$ and a rational map  $$ F : X \times T
\longrightarrow \bold P(\Cal Q)  $$ To construct $\Cal Q$ we fix a Poincar\'e
bundle $$\Cal P \rightarrow Pic^{2g+2}(C) \times C$$ and consider the vector
bundles $$ \Cal F = p_*\Cal P \quad \text {and } \quad \Cal G = p_*(\Cal
P^{\otimes 2}) \tag 2.21 $$ where $p$ is the first projection. Then $ \bold
P(\Cal F)= C^{(2g+2)} $, the fibre of $\Cal F$ is $H^0(L)$ and the fibre of
$\Cal G$ is
$H^0(L^{\otimes 2})$ (at the point $L \in Pic^{2g+2}(C)$). Moreover there
exists a uniquely defined
morphism $$ \nu : Sym^2 \Cal F \longrightarrow \Cal G \tag 2.22  $$\noindent
such that   $\nu_L:
(Sym^2 \Cal F)_L \rightarrow \Cal G_L$ is the usual multiplication map $$\nu_L:
Sym^2 H^0(L)
\longrightarrow H^0(L^{ \otimes 2}) $$   $\nu$  Since $deg(L) \geq 2g+1$
$\nu_L$
is always surjective, hence
  $$ \Cal V = Ker \nu \tag 2.23 $$\noindent is a vector bundle and one computes
$$
\text { rank of $\Cal V$ =   $\binom {g+1}2 +1 $}
$$
\noindent We define now a squaring map
$$ \gamma: T \longrightarrow Pic^{2g+2}(C)  $$
\noindent by setting
$$
 \quad \gamma (t) = \omega_C(2t)
$$
\noindent then we take the pull-back of $\Cal V$ by $\gamma$ and give the
following
\proclaim {(2.24)DEFINITION} $\Cal Q = \gamma^* Ker \nu $ is the
fundamental bundle. \endproclaim
\noindent Considering the multiplication map
$$
\mu_t : Sym^2 H^0(\omega_C(2t)) \longrightarrow H^0(\omega_C^{\otimes 2}(4t))
\tag 2.25
$$
\noindent it is clear that $Ker \mu_t$ is the fibre of $\Cal Q$ at $t$,
therefore
$$
\Cal Q_t = Ker \mu_t = H^0(\Cal I_t(2)) \tag 2.26
$$
\noindent Now we want to construct $F$: we fix the open set $$ U =  U_j \cap
(X-Sing X) \tag2.27 $$
\noindent where $U_j$ is defined as in the previous section. Let $(z,t) \in U$
then $z$ is the moduli
point of a stable $\xi$ such that $h^0(\xi(t)) = 4$ (prop.2.12). Putting as
usual $\Cal E$ =
$\xi(t)$, $V = H^0(\xi(t))$   we construct from $(\Cal E,V)$ the quadratic form
$q(\Cal E,V) \in
H^0(\Cal I_t(2))$. Since $\xi$ is stable and $(z,t) \in U_j$ there is no
subline bundle $L$ of $\Cal
E$ with $h^0(L) \geq 2$ (prop. 2.15), in particular $q(\Cal E,V)$ is not zero
(prop.1.11).
Therefore its zero locus  $Q(\Cal E,V)$ is a quadric containing $C_t$ i.e. a
point of $\bold P \Cal Q_t$. $Q(\Cal E,V)$ depends only on $(z,t)$, we set
$$ F(z,t) = Q(\Cal E,V)  $$\noindent this defines a rational map  $$ F: X
\times T
\longrightarrow \bold P \Cal Q $$\noindent which is regular on $U$.  \proclaim
{(2.28)DEFINITION}ÊWe will say that $F: X \times T \longrightarrow \bold P \Cal
Q$ is the fundamental map. The restriction of $F$ to $X \times t$ will be
denoted by $$ F_t: X  \longrightarrow \bold P \Cal Q_t
$$
\endproclaim
\bigskip \noindent
{\bf 3. Geometry of the fundamental map: rank 6 quadrics \rm } \bigskip
\noindent The family
of maps $F_t: X \rightarrow \bold P \Cal Q_t$ has the following important
feature: let $$\theta: X \rightarrow \bold PH^0(\Cal L)^*$$\noindent be the
morphism
associated to the generalized theta divisor $\Cal L$ then $F_t$ is defined by a
linear subsystem of $\mid \Cal L \mid$. In other words $\lbrace F_t(X), \quad t
\in T
\rbrace $ is a complete family of linear projections of $\theta(X)$. This is
shown in section 4 and it will be used to discuss the very ampleness of
$\Cal L$.  Now we want to study the map  $F_t$: this provides a natural
description of $X$ in terms of rank 6 quadrics. Let $$ U = U_j \cap X-Sing X
\quad U_t = U \cap X
\times t \tag 3.1  $$\noindent  by proposition 2.11(i) $U_t$ is not empty for
all $t\in T$. $F_t$
is regular on $U_t$, we consider  the closure   $$ Z_t \subset \bold P\Cal Q_t
\tag 3.2$$\noindent  of
 $F_t(X)$ and the scheme theoretic intersection
 $$
W^{(r)}_t = {\Cal P}^{(r)}_t \cdot \bold P(\Cal Q_t) $$\noindent
where $\Cal P_t^{(r)} \subset \bold PH^0(\Cal O_{\bold P^{g+2}_t}(2))  $ is the
variety  of
all quadrics in $\bold P^{g+2}_t$ having rank $\leq r$. It is clear that
$$
Z_t \subseteq W^{(6)}_t
$$
\noindent moreover the expected dimension of $W^{(6)}_t$ is $3g-3$ and this is
also the dimension
of $X$: the situation is described by the next \proclaim {(3.3)THEOREM} Assume
$g \geq 3$, then \par
\noindent (1) the general point of $Z_t$ is a quadric of rank 6 \par \noindent
(2) the rational map
$F_t : X \longrightarrow Z_t$ has degree two \par \noindent (3) $Z_t$ is a
reduced irreducible
component of $W_t$. \endproclaim \demo {Proof} Preliminarily we show that the
map $F_t: X \rightarrow
Z_t$ has degree $\leq 2 $:  Let $z \in U_t$, then $z$ is the moduli point of a
stable bundle
$\xi$; since $z \in U_t$ $\xi(t)$ is globally generated, there is no subline
bundle $L$ of $\xi(t)$ with $h^0(\xi(t)$ and $h^0(\xi(t) = 4$(prop. 2.12,
2.15). Therefore (by prop.
1.11)
   $$ \text { $ 5 \leq rank F_t(z) \leq 6$}Ê\quad \text {   $ Sing F_t(z) \cap
C_t =
\emptyset $} $$\noindent  Since these conditions hold  we can apply to $F_t(z)$
proposition 1.19 which says that there exist at most two isomorphism classes of
pairs
$(\Cal E,V)$ such that $Q(\Cal E,V)$ = $F_t(z)$. Hence $deg F_t \leq 2$ and
moreover  $$ dim Z_t =
3g-3 $$ because $F_t$ is of finite degree. We can now prove our statements:
assume (1) does not hold,
then every point of $Z_t$ is a quadric of rank $\leq 5$ and   $$ Z_t \subseteq
W_t^{(5)}  $$\noindent
Fix $Q$ = $F_t(z)$ with $z \in U_t$ then $Q$ has rank exactly 5.  Consider  the
projectivized tangent
space $T^{(5)}_Q$  to  $W_t^{(5)}$ at the point $Q$, it is standard that
$$T^{(5)}_Q = \bold PH^0(\Cal I_{Sing Q\cup C}(2))$$\noindent where $\Cal
I_{Sing Q \cup C}(2))$ is
the ideal sheaf of $Sing Q \cup C$. Since $Z_t$ is a closed  subset of
$W^{(5)}_t$ it follows $$ dim
T^{(5)}_Q \geq 3g-3 = dim Z_t $$ Let us show that this is impossible  for $g
\geq 3$:since $Sing Q
\cap C$ is empty  the general maximal linear subspace $\Lambda$ of $Q$ does not
intersect $C$. Then,
by proposition 1.26, $$ h^0(\Cal I_{\Lambda \cup C}(2)) = h^0(\Cal E) - 3
$$\noindent where $\Cal
I_{\Lambda \cup C}$ is the ideal sheaf of $\Lambda \cup C$ and $(\Cal E,V)$ is
a pair defining $Q$.
Since $Q$ = $F(z,t)$ with $z \in U_t$ we can assume $\Cal E = \xi(t)$ where
$\xi$ is the bundle
considered at the beginning of this proof. But then $h^0(\Cal E) = 4$ so that
$$ h^0(\Cal I_{\Lambda
\cup C}(2)) = 1 $$\noindent On the other hand we can consider the exact
sequence $$ 0 \rightarrow
H^0(\Cal I_{\Lambda \cup C}(2)) \rightarrow H^0(\Cal I_{Sing Q \cup C}(2))
\longrightarrow H^0(\Cal
J(2))  $$\noindent where  $\Cal J$ is the ideal of $Sing Q$ in $\Lambda$. We
have $h^0(\Cal I_{Sing Q
\cup C}(2)) \geq 3g-2 = dim T^{(5)}_Q + 1$, from  $dim Sing Q = g-3$ and $dim
\Lambda = g-1$ we
compute immediately  $h^0(\Cal J(2)) = 2g-1$. Thus we  finally obtain $$
h^0(\Cal
I_{\Lambda \cup C}(2)) \geq g-1 $$\noindent
which is impossible for $g \geq 3$. This shows (1), to show (2)
observe that, if $F_t(z)$ has rank 6, there exist exactly two isomorphism
classes of pairs $(\Cal
E,V)$ defining $Q$. To show (3) let
$Q$ = $F_t(z)$ with $z \in U_t$ and rank of $Q$ = 6. We consider  the
projectivized tangent space to $W^{(6)}_t$ at $Q$  that is $$ T^{(6)}_Q =
\bold P H^0(\Cal I_{Sing Q \cup C}(2)) $$\noindent Then we use the previous
sequence 3.14: this time
we have $dim Sing Q = g-4$, $dim \Lambda = g-1$ so that $h^0(\Cal J(2)) =
3g-1$. On the other hand,
choosing $\Lambda$ general in at least one of the two rulings of maximal linear
subspaces, we have
$h^0(\Cal I_{\Lambda \cup C}(2)) = 1$ as above. Therefore $h^0(\Cal I_{Sing Q
\cup C}(2)) \leq 3g-2$
and $dim T^{(6)}_Q \leq 3g-3$. Since $dim T^{(6)}_Q \geq dim Z_t =3g-3 $ it
follows $dim T^{(6)}_Q =
3g-3$. Hence $Z_t$ is an irreducible component of $W^{(6)}_t$ which is smooth
at $Q$.
 \enddemo
Let $j: X \times T \rightarrow X \times T$ be the fundamental involution, it
follows from the
previous proof that generically
$$ F^{-1}(F(z,t)) =
\lbrace (z,t) \quad , \quad j(z,t) \rbrace $$ \noindent
hence we have also obtained:
 \proclaim
{(3.5)THEOREM} Assume $gÊ\geq 3$, then:
\par \noindent (1) the fundamental map $F$ has degree two; \par \noindent (2)
the birational
involution induced by $F$ is the fundamental involution  $j$. \endproclaim
Assume $g \geq 3$, then $Z_t$ is a family of quadrics of even rank 6:  a  well
known  construction
defines a finite double covering $$
\tilde{F_t}: {\tilde Z}_t \longrightarrow Z_t \tag 3.6
$$
\noindent where, as a set, ${\tilde Z}_t$ is the family of pairs $(Q,\Lambda)$
such that $Q \in
Z_t$,  $\Lambda$ is a connected component of the variety of linear subspaces of
$Q$ having
codimension 2 and $\tilde {F_t}(Q,\Lambda) = Q$. In particular $\tilde{F_t}$ is
\'etale at
each $Q \in Z_t$ having rank 6.  \proclaim {(3.7)THEOREM} $F_t = \tilde {F_t}
\cdot \pi_t$
where $\pi_t: X \longrightarrow \tilde {Z_t} $ is a birational morphism
\endproclaim \demo
{Proof} $F_t(z)$ = $Q(\Cal E,V)$ for a given pair $(\Cal E,V)$. As in section 1
we consider the
grassmannian $G_V^*$ and the universal bundle $\Cal U_V$: one of the two
families of planes in
$G_V^*$ is the family of the zero sets of the non zero global sections of $\Cal
U_V^*$.  Recall that
$Q(\Cal E,V)$ = ${\delta_V}^{-1}(G_V^*)$ where $\delta_V$ is the linear
projection defined in 1.5.
Taking the inverse images by $\delta_V$ of the planes of this family we obtain
a
connected component $\Lambda$ of the variety of codimension 2 linear subspaces
of
$Q(\Cal E,V)$. One defines $\pi_t$ by setting $\pi_t(z,t) =
(Q,\Lambda)$. \enddemo We finish this section by an auxiliary result to be used
later: \proclaim {(3.8)PROPOSITION} Let $\xi$ be a stable bundle, $z \in X$ its
moduli point, $t \in T$. Assume $h^0(\xi(t) = 4$, $\xi(t)$ globally generated
and
$F(z,t)$ of rank 6. Then $F(t,z)$ is a smooth point of $Z_t$ and the tangent
map $dF$
is injective at $(z,t)$ as well as $(dF_t)$ at $z$. \endproclaim \demo {Proof}
Let $Q $ =
$F(z,t)$: in theorem 3.3 we proved  $dim Z_t$ = $3g-3$ exactly by showing that
the tangent space
$T_{Z_t,Q}$ has dimension $3g-3$ (under the assumptions we have). Hence $Q$ is
non singular for
$Z_t$. Since $\xi(t)$ is globally generated and $Q$ has rank $\geq 5$ it
follows from our usual
argument (prop. 1.19) that $z$ is a connected component of $F^{-1}(z)$. $z$ is
smooth because $\xi$
is stable, since $F_t = \tilde {F_t} \cdot \pi_t$ and $\tilde {F_t}$ is \'etale
over $Q$ because the
rank of $Q$ is 6 we conclude that $(dF_t)_z$ is injective. Finally let  $p:
\bold PÊ\Cal Q
\rightarrow T$ be the natural map, since $p \cdot F = id_T$ then $dp_Q \cdot
dF_{(z,t)}:
T_{X,z}\oplus T_{T,t} \rightarrow T_{T,t}$ is the obvious projection, on the
other hand
$dF_{(z,t)}/T_{X,z} = (dF_t)_z$. This implies that  $dF_{(z,t)}$   is
injective.      \enddemo
\bigskip \proclaim {\bf 4. Geometry of the fundamental map:  the generalized
theta
divisor \rm} \endproclaim \bigskip \noindent In this section, and throughout
all
the paper, we fix the following notations \proclaim {(4.1) } \par \noindent
$C^{(n)}$ = n-th symmetric product of the curve $C$, \par \noindent
$J=Pic^0(C)$,\par \noindent  $\Theta$ = a symmetric theta divisor in $J$, \par
\noindent $X = SU(2,C)$,\par \noindent $\Cal L$ = generalized theta divisor of
$SU(2,C)$ \par \noindent Moreover we will set
 $$ S = C^{(2)} $$ \noindent
 and
 $$ T = Pic^2(C) $$ \endproclaim \par \noindent
The vector space $H^0(\Cal L)^*$ will be canonically identified to $H^0(\Cal
O_J(2\Theta))$ as in [B1].
Though it is actually not needed we will assume in this section that $C$
is \it not hyperelliptic: \rm  this  simplifies somehow the exposition. Our
first purpose is to construct  a suitable map of
vector bundles over $T$ $$ \lambda : H^0(\Cal L)^* \otimes \Cal O_T
\longrightarrow \Cal Q $$\noindent where $\Cal Q$ is the \it fundamental bundle
\rm defined in section
3. In order to do this fix $t \in T$ and consider the morphism
 $$ \alpha_t : S \longrightarrow J  $$\noindent such that $$ \alpha_t(x+y) =
t-x-y$$\noindent
$\alpha_t$ is just the Abel map multiplied by -1, its image in $J$ is $t-S$; we
denote it by $$  S_t  \tag 4.2 $$\noindent This defines the correspondence $$
\Cal S =
\lbrace (t,e) \in T \times J | e \in  S_t \rbrace \tag 4.3 $$ together with its
two natural projections $$ \CD T @<{\pi_1}<< {\Cal S} @>{\pi_2}>> J \\ \endCD
$$\noindent
Applying $\pi_{1*} \pi_2^*$ to the evaluation map $e :\Cal O_J \otimes H^0(\Cal
O_J(2\Theta))  \longrightarrow \Cal O_J(2\Theta)$
we obtain the morphism of sheaves
$$
\pi_{1*}\pi_2^*(e): \pi_{1*}\pi_2^* \Cal O_J \otimes H^0(\Cal
O_J(2\Theta)) \longrightarrow \pi_{1*}\pi_2^* \Cal O_J(2\Theta)
$$
\noindent let us set
$$
\hat {\Cal Q} = \pi_{1*}\pi_2^* \Cal O_J(2\Theta) \tag 4.4
$$
\noindent and
$$
\hat{\lambda} = \pi_{1*}Ê\pi_{2}^*(e) \tag 4.5
$$
\noindent since  $\pi_{1*}\pi_2^* \Cal O_J \otimes H^0(\Cal O_J(2\Theta)) =
\Cal O_T \otimes
H^0(\Cal O_J(2\Theta))$ we have defined  a morphism  $$
\hat{\lambda} = \pi_{1*}\pi_2^*(e) : H^0(\Cal O_J(2\Theta)) \otimes \Cal O_T
\longrightarrow \hat {\Cal Q}
 $$
\noindent Note that $\pi_1^* t =  S_t$ so that  the stalk of $\hat {\Cal Q}$ at
$t$ is  $$
\hat {\Cal Q}_t = H^0(\Cal O_{S_t}(2\Theta))
$$
\noindent and the map  $$\hat{\lambda}_t : H^0(\Cal O_J(2\Theta)) \rightarrow
\hat
{\Cal Q}_t $$\noindent is just the restriction  $$ H^0(\Cal O_J(2\Theta ))
\longrightarrow H^0(\Cal O_{S_t}(2\Theta)) $$\noindent
\proclaim {(4.6)DEFINITION} We will say that $\hat {\lambda} $ is the
restriction map. \endproclaim \noindent If $dim (\hat{\Cal Q}_t)$ is a constant
function $\hat{\Cal Q}$ is a vector bundle: this is ensured by the next
proposition 4.9.

\proclaim {(4.7)LEMMA } Let $d = \Sigma z_i \in Div(C) $, $D=
\Sigma (z_i+C)$ be the corresponding divisor in $S$, $\Delta$ the diagonal in
$S$. Then: \par \noindent
(1) $\Cal O_{\Delta}(D) \cong \Cal O_C(2d)$ \par \noindent
(2) there exists a canonical isomorphism  $\psi: Sym^2H^0(\Cal
O_C(d)) \rightarrow H^0(\Cal O_S(D))$ such that\par \noindent (3)
composing the restriction $
\rho : H^0(\Cal O_{S}(D) \longrightarrow H^0(\Cal O_C(2d))
$
with $\psi$ we  obtain the multiplication map
$$
\mu : Sym^2H^0(\Cal O_C(d)) \longrightarrow H^0(\Cal O_C(2d))
$$
\noindent In particular $\rho$ is surjective if $deg(d) \geq 2g+1 $.
\endproclaim
\demo {Proof} It is elementary to check that  $(z+C)$ and $\Delta$ intersect at
the unique point $2z \in S$ with multiplicity  2: this implies (1). To show (2)
consider the natural involution $i: C \times C \rightarrow C\times C$
together with the quotient map $\pi: C \times C \rightarrow S $ and the
projections $p_i: C\times C \rightarrow C$. From  $p_1^*\Cal O_C(d) \otimes
p_2^* \Cal O_C(d) $ $\cong $ $\pi^* \Cal O_S(D)$ we obtain the isomorphism
$$ Ê\tau: H^0(\Cal O_C(d)) \otimes H^0(\Cal
O_C(d)) \rightarrow H^0(\pi^*\Cal O_S(D))  $$\noindent  such that $ \tau (s_i
\otimes s_j)  =
p_1^*s_i \otimes p_2^* s_j $. On the other hand we have the involution
$\iota^*$ on $
H^0(\pi^*\Cal O_S(D))$ defined by $\iota$. Let $ j = ({\tau}^{-1} \cdot
{\iota}^* \cdot {\tau})$:
observe that $j$ is  the natural involution on $H^0(\Cal O_C(d))Ê\otimes
H^0(\Cal O_C(d))$ sending
$s_i \otimes s_j$  to $s_j \otimes s_i $. \par \noindent Considering the +1
eigenspaces we have
$Sym^2H^0(\Cal O_C(d))$ for $j$ and  $\pi^* H^0(\Cal O_S(D))$ for $\iota^*$.
Restricting
$\tau$ to $Sym^2 H^0(\omega_C(2t))$ we obtain the diagram
  $$ \CD { Sym^2 H^0(\Cal O_C(2))} @>{\tau}>> {\pi^*H^0(\Cal
O_{S}(D))} @<{\pi^*}<< {H^0(\Cal O_S(D))} \\
\endCD $$
\noindent from this (composing properly the maps) we obtain the isomorphism
$\psi$ required in
(2). The equality $\mu = \rho \cdot \psi$ can be  checked on  vectors $s_i
\otimes s_j + s_j \otimes
s_i$. Finally $\mu $ is surjective if $deg d \geq 2g+1$ so that the same holds
for $\rho$.
\enddemo \noindent We denote by $$ \Theta_t \in Div(S) $$\noindent a divisor
such that $$
\Cal O_S(\Theta_t) \cong \alpha_t^*\Cal O_J(\Theta) $$\noindent for brevity we
omit the proof
of the following \proclaim {(4.8)LEMMA} $\Cal O_S(2\Theta_t) \cong \Cal
O_S(\Sigma
(z_i+C) - \Delta) $ where $\Cal O_C(\sum z_i) \cong \omega_C(2t)$. \endproclaim
 \proclaim {(4.9)PROPOSITION} $\forall t
\in T$ we have \par \noindent (1)$h^0(\Cal O_{S_t}(2\Theta))$ = $\binom
{g+1}2 +1$  \par \noindent (2) $h^1(\Cal O_{S_t}(2\Theta))=0$
\endproclaim \demo {Proof} (1) By the previous two lemmas we have the exact
sequence $$ 0
\rightarrow \Cal O_{S}(2\Theta_t) \rightarrow \Cal O_{S}(\Sigma (z_i +C))
\rightarrow \Cal O_C(2\Sigma z_i) \rightarrow 0 \tag 4.10 $$\noindent with
$\Cal
O_C(\Sigma z_i) \cong \omega_C(2t)$. Consider the restriction map  $$ \rho
:H^0(\Cal O_{S}(\Sigma (z_i+C)) \longrightarrow H^0(\Cal O_C(2\Sigma z_i))
$$\noindent up
to the  isomorphism $\psi$ given in lemma  4.7  $\rho$ is just the
multiplication map $$ \mu_t :Sym^2 H^0(\omega_C(2t)) \longrightarrow
H^0(\omega_C^{\otimes 2}(4t)) \quad \tag 4.11 $$\noindent
since $deg(\Sigma z_i) \geq 2g+1$ $\mu_t$ is surjective and  $Ker \mu_t \cong
H^0(\Cal O_S(2\Theta_t))$. Then we can compute $h^0(\Cal O_S(2\Theta_t))$ =
$\binom {g+1}2 +1$. \par \noindent (2) By Riemann-Roch $h^0(\Cal
O_{S}(2\Theta_t))= \binom {g+1}2 +1 + h^1(\Cal O_{S}(2\Theta_t)) - h^2(\Cal
O_{S}(2\Theta_t)$. By Serre duality $h^2(\Cal O_{S}(2\Theta_t))=h^0(\Cal
O_{S}(K_{S}-2\Theta_t))$.
 It is easy to see that $h^0(\Cal O_{S}(K_{S} - 2\Theta_t)) = 0$;
indeed $z+C$ is nef and $(z+C)$ $\cdot$ $ (K_{S}-2\Theta_t)$ = $-3$. Hence, by
(1), $h^2(\Cal O_{S}(2\Theta_t))= h^1(\Cal O_{S}(2\Theta_t))=0$.\enddemo
\proclaim {(4.12)COROLLARY} $\hat {\Cal Q}$ is a vector bundle of rank $\binom
{g+1}2 + 1$ over $Pic^2(C)$ \endproclaim \noindent
Note that this is also the rank of the \it fundamental bundle \rm $\Cal Q$
constructed in
section 3. Moreover we have the exact sequence
$$
\CD
0 @>>> {H^0(\Cal O_S(2\Theta_t))} @>{\sigma_t}>> {Sym^2 H^0(\omega_C(2t))}
@>{\mu_t}>>
{ H^0(\omega_C^{\otimes 2}(4t))} @>>> 0 \\
\endCD
$$\noindent
where $\mu_t$ is the multiplication map: this is obtained from  the sequence
4.10 applying
lemma 4.7. Since the fibre of $\Cal Q$ on $t$ is $\Cal Q_t$ = $Ker \mu_t$ we
have a natural
isomorphism
 $$ \sigma_t : \hat{\Cal Q}_t \longrightarrow \Cal Q_t \tag 4.13  $$
 \proclaim {(4.14)PROPOSITION} Up to tensoring by a line bundle there exists an
isomorphism $ \sigma :
\hat{\Cal Q} \longrightarrow \Cal Q $ having as its fibrewise maps the previous
isomorphisms
$\sigma_t$ ($\forall t \in T$). \endproclaim \noindent Since this is
essentially not needed we omit
the proof. Composing the restriction map $\hat{\lambda}$ with $\sigma$ we
obtain the morphism of
vector bundles  $$\lambda = \sigma \cdot \hat{\lambda}: H^0(\Cal O_J(2\Theta))
\otimes \Cal O_T
\rightarrow \Cal Q$$\noindent  From $\lambda$ we construct  the diagram $$ \CD
{X \times T} @>{\theta
\times id}>>{\mid 2\Theta \mid \times T } @>{\overline{\lambda}}>> {\bold P
\Cal Q} \\ \endCD
$$\noindent where  $\theta: X \rightarrow \mid 2\Theta \mid$ is the morphism
associated to the
generalized theta divisor and $\overline{\lambda}$ is the projectivization of
$\lambda$.\par
\noindent This defines the rational map $$ \Phi = \overline {\lambda} \cdot
(\theta \times id_T)
X \times T \longrightarrow \bold P \Cal Q \tag 4.15 $$\noindent on the other
hand we have from section
2 the fundamental map $$ F: X \times T \longrightarrow \bold P \Cal Q \tag
4.16$$ \noindent The main
step is now the following
  \proclaim {(4.17)THEOREM} $F$ and $\Phi$ are the same map. \endproclaim
 In order to show the theorem we need some preparation. \proclaim
{(4.18)LEMMA} Let $C \subseteq \bold P^r$ be a not degenerate smooth curve,
$Sec(C)$ the variety of its bisecant lines. Then there  is no quadric
containing $Sec(C)$. \endproclaim \demo {Proof} Assume there exists a quadric
$Q$ containing $Sec(C)$ then $Q$ contains the linear span of $<d>$ for every
$d \in C^{(3)}$ (indeed $<d>$ is a 3-secant plane or a 3-secant line to $C$: in
the former case $Q$ contains three lines of $<d>$ hence $<d>$). Using the same
argument one can show by induction on $k \geq 3$ that $Q$ contains the linear
span of $<d>$ for
every $d \in C^{(k)}$. Since $C$ spans $\bold P^r$ this implies $Q$ = $\bold
P^r$: obviously a contradiction. \enddemo \noindent Let
$$
Sec(C_t) \subset \bold P^{g+2}_t
 \tag 4.19 $$\noindent be the variety of bisecant lines to $C_t$: $\forall x,y
\in C$ we
denote by $$ \overline {xy}_t $$\noindent the bisecant line to $C_t$ joining
$x$ to $y$. We consider
the incidence correspondence $$ \Sigma_t = \lbrace (x+y,z) \in S \times \bold
P^{g+2}_t / z \in
\overline {xy}_t \rbrace $$\noindent together with its two projections $$  \CD
S @< {\alpha}<< {\Sigma_t} @>{\beta}>>{Sec (C_t)} \\
\endCD \tag 4.20
$$
\noindent Let
$$
E = \lbrace (x+y,z) \in \Sigma_t / \quad \text { $x = z$ or $y = z$} \rbrace
\tag 4.21 $$
\noindent $E$ is a copy of $C \times C$ and it is contracted to $C_t$ by
$\beta$.
\proclaim {(4.22) LEMMA}\par \noindent (1) $\beta^* \Cal O_{Sec (C_t)}(2)
\otimes \Cal O_{\Sigma_t}(-E) \cong \alpha^*\Cal O_S(2\Theta_t) $ \par
\noindent
(2)Let  $s \in
H^0(\Cal O_S(2\Theta_t))$, $s \neq 0$  then
  $$ \alpha^* div(s) - E = \beta^*Q $$\noindent
where $Q$ is the quadric defined by $\sigma_t(s)$. \endproclaim
\demo {Proof} Let $s \in \hat {\Cal Q}_t = H^0(\Cal O_S(2\Theta_t))$, $s \neq
0$: at first we construct
  $$ q = \sigma_t(s)
\tag 4.23 $$ \noindent (up to a non zero constant factor). For this it suffices
to use
the proof of lemma 4.7 and to recall the definition of $\sigma_t$ as it has
been
given in 4.13:  the  sequence 4.10 defines an injection $H^0(\Cal
O_S(2\Theta_t)) \rightarrow H^0(\Cal O_S(D))$. Let
 $$ r $$\noindent be the image of $s$ in $H^0(\Cal O_S(D))$: using the quotient
map $\pi
: C \times C \rightarrow S$ we lift $r$ to a global
section  $$ b = \pi^*r \in H^0(\pi^*\Cal O_S(D)) $$\noindent
 As in the proof of lemma 4.7 we have the standard
identifications $$ H^0(\omega_C(2t)) \otimes H^0(\omega_C(dt)) = H^0(\pi^*\Cal
O_S(D)) $$ and $$ Sym^2H^0(\omega_C(2t)) = \pi^*H^0(\Cal O_S(D))
$$
\noindent Hence we can view $b$ as a symmetric bilinear form on
$H^0(\omega_C(2t))^*$: let $q$ be its associated quadratic form then
 $$ q = \sigma_t(s)$$\noindent To prove (1) and (2),  it suffices
to show that $$ div(\beta^* q_s) - E = \alpha^* div(s) $$
Let us check this equality set theoretically: by definition of $b$ we have $$
x+y \in div(r) \Longleftrightarrow <x,y> \text {is an isotropic space for
$b$} \tag 4.24 $$\noindent (where $<x,y> $ $\subset$ $ H^0(\omega_C(2t))^*$
denotes
the vector space having as its projectivization the linear span of $x,y$ in
$\bold P^{g+2}_t$). Furthermore $div (r) = div(s) + \Delta$ so that  $b(x,x)
= 0$ $\forall x \in C$ and, as it must be, $q$ vanishes on $C_t$; let $Q$ be
the quadric defined by $q$ and $$
R = Q \cap Sec(C_t)
$$\noindent
since $Q$ contains $C_t$ one can easily show that $R$ is union of bisecant
lines
$\overline {xy}_t$. Assume $x \neq y$ then, using the previous equivalence
and $div(r) = div(s) + \Delta$, it follows   $$ \overline{xy}_t \subset Q
\Longleftrightarrow b_s(<x>,<y>) \Longleftrightarrow x+y \in div(s) $$\noindent
Hence
$$
\beta^{-1}(R) = \alpha^{-1}(div(s)) \cup E \tag 4.25
$$
\noindent For a general $s$ $div(s)$ is reduced, hence $\alpha^* div(s) =
alpha^{-1}(div(s)$; moreover a standard computation in $Num \Sigma_t$ shows
that $\alpha^*div(s)$ is numerically equivalent to $\beta^*q - E$. From
this and the set theoretical equality the proof follows. \enddemo
 \noindent Now we can give a \proclaim{(4.26) PROOF OF THEOREM 5.20}
\endproclaim \demo {Proof} Let us recall that: \par \noindent (1)  $F$ is the
composition of the following maps
$$
\CD
{X \times T} @>{\theta \times id_T}>> {\mid 2\Theta \mid \times \times T}
@>{\hat {\lambda}}>> {\bold P \hat {\Cal Q}} @>{\sigma}>> {\bold P \Cal Q} \\
\endCD
$$
\noindent (2) given $z \in X$ we have
$$\theta(z) = \Theta_{\xi} = \lbrace e \in J / h^0(\xi(e)) \geq 1 \rbrace $$
where $\xi$ is any bundle with moduli point $z$. \par \noindent Assume $(z,t)$
is in the domain of $F$, since $\hat{\lambda}_t$ is  the
restriction map $H^0(\Cal O_J(2\Theta)) \rightarrow H^0(\Cal O_{S_t}(2\Theta))$
it follows that $$ \hat {\lambda}(\theta(z),t)$$ is the curve
   $$ \Theta_{t,z} = \theta(z) \cdot S_t $$
\noindent let $\Theta_{z,t} = div(s)$, $q = \sigma_t(s)$, then $$ Q = F(z,t)
\in \bold
P \Cal Q_t $$ is the quadric defined by $q$. In particular it follows that:
\par \noindent $(z,t)$ is in the domain of $F$ $\Longleftrightarrow$
$\Theta_{t,z}$ is well defined as a divisor in $S_t$ $\Longleftrightarrow$
$\theta(z)$ does not contain $S_t$ \par \noindent
Actually it is possible to show that these conditions are equivalent to
$h^0(\Cal E) = 4$ and $q(\Cal E,V)$  not identically zero, where $ \Cal E =
\xi(t)$, $ V = H^0(\Cal E) $; in any case, for a general
 $(z,t)$ in the domain of
$F$, the latter conditions are satisfied by $\xi$. This gives to us a second
quadric $$ Q' = Q(\Cal E,V)$$\noindent which is defined defined by the pair
$(\Cal E,V)$
as in section 1. We must show $Q$ = $Q'$. Since no quadric contains  $Sec
(C_t)$ it
suffices to show $Q \cdot Sec(C_t) = Q'\cdot Sec(C_t)$ that is
 $$ \beta^*Q-E = \beta^*Q'-E
$$\noindent  At first we check this equality set theoretically: since $Q$ is
defined by
$\sigma_t(s)$ with $div (s) = \Theta_{t,z}$ we have from the previous lemma
$$ \alpha^*\Theta_{t,z}-E = \beta^*Q $$
\noindent which means
$$
x+y \in \Theta_{t,z} \Longleftrightarrow \overline{xy}_t \subset Q $$
\noindent on the other hand
$$
x+y \in \Theta_{z,t} \Longleftrightarrow x+y \in \Theta_{\xi} \cap S_t
\Longleftrightarrow h^0(\xi(t-x-y)) \geq 1
$$
\noindent by definition of $\Theta_{\xi}$. Now recall that $Q' = Q(\Cal E,V)$
so that,
using the same notations of section 1, we have the linear projection
$$
\delta_V : Q_1 \longrightarrow G_V^*
$$
\noindent such that $\delta_V/C_t$ is the Gauss map $g_V$. It is a standard
exercise to
check that
$$
h^0(\xi(t-x-y)) \geq 1 \Longleftrightarrow \delta_V(\overline{xy}_t) \subset
G_V^* \Longleftrightarrow \overline{xy}_t subset Q'
$$
\noindent Therefore $Q \cap Sec(C_t) = Q' \cap Sec(C_t)$. To complete the proof
 it suffices to show that $Q \cap Sec(C_t)$
= $Q \cdot Sec(C_t)$ which is true if  $\Theta_{z,t}$ is reduced. Now, for
a general $z$, $\Theta_{\xi}$ is reduced; moreover, by transversality of
general
translate, $\Theta_{\xi} \cdot S_t$ is reduced too for a general $t$.
\enddemo
\bigskip
\proclaim { \bf 5. Geometry of the fundamental involution: technical
lemmas \rm} \endproclaim \noindent Let $j: X \times T \rightarrow X \times
T$ be the fundamental involution: from now on we will denote respectively by
$$
I_j \text { and $B_j$ }\tag 5.1
$$
the indeterminacy locus of $j$ and the maximal (with respect to inclusion) open
subset
along which $j$ is biregular. The proof of  our main theorem relies on a weaker
version
of rigidity lemma which is given in the next section. This rigidity argument
can be
certainly applied if: \it for every $o \in X - Sing X$ \rm \bigskip \noindent
(1)   $  o  \times T \cap I_j$ has codimension
$\geq 2$ in $T$ \par \noindent (2)  $ o  \times T \cap B_j$ is not
empty. \bigskip \noindent   Since we did not come to a complete proof of
this statement, we found more convenient showing a weaker result which is
sufficient for our
purposes: this is the content of the next theorems 5.6, 5.7,5.8. Let $\xi$ be a
semistable rank 2
vector bundle, $ det \xi = \omega_C$, we define from $\xi$ the  following
subsets of $T$ \bigskip
\noindent (5.2) \par \noindent
(1) $V_{00}(\xi) = \lbrace t \in T/ \quad h^0(\xi(t-x-y) \geq 1$ \it for all
\rm $x$, $y$
$\in C \rbrace$, \bigskip \noindent
(2) $V_0(\xi) = \lbrace t \in T/ \quad h^0(\xi(t-2x) \geq 1$ \it for all \rm
$x$,
$\in C \rbrace$ \bigskip \noindent
(3) $V_m(\xi) = \lbrace t \in T/\quad  h^0(\xi(t-mx)) \geq
4-m$ for some $x \in C \rbrace $, $m = 1,2,3$ \bigskip \noindent
(3) $W(\xi) = \lbrace t \in T/ \quad h^0(\xi(t-L)) \geq 1$ for some $L \in
W^1_d$,
$d\leq g \rbrace $ \bigskip \noindent  where $W^1_d \subset Pic^d(C) $ is the
\it
Brill-Noether locus \rm of line bundles $L$ having degree $d$ and $h^0(L) \geq
2$.
In addition, considering the following subset of $V_1(\xi) \cap V_3(\xi)$ will
be quite important:
$$
V_{13}(\xi) = \lbrace t \in T/ \quad h^0(\xi(t-x))Ê\geq 3 ,\quad h^0(\xi(t-3x))
\geq 1 \text { for \it the same \rm $x \in C$ }\rbrace \subset V_1(\xi) \cap
V_3(\xi) \tag 5.3 $$By semicontinuity these sets are closed in $T$, it is easy
to see that they
depend only on the moduli point $[\xi]$ of $\xi$; let
$$
I(o) = V_0(\xi) \cup V_{13}(\xi) \cup W(\xi)   \tag 5.4
$$
if $[\xi]Ê= o$ and let
$$
I(o_1 \dots o_r) = I(o_1) \cup \dots \cup I(o_r)
$$
\noindent for points $o_1 \dots o_r \in X$
\proclaim {(5.5) DEFINITION}Ê$I(o_1, \dots , o_r)$ is the special set of
$\lbrace
o_1 \dots o_rÊ\rbrace$. \endproclaim
\proclaim {(5.6) THEOREM} Let $o, \overline o \in X $, assume: (\bf 1 \rm)
$ o \times T \cap B_j $  is not empty; (\bf 2 \rm)$ j( o \times T) \subset
\overline o  \times T $.
Then    $$
I_j \cap o \times T \subseteq  o  \times I(o,\overline o)
$$
\noindent and the same holds for $\overline o$.
\endproclaim
\proclaim {(5.7)THEOREM} Let $o = [\xi] \in X$ then $ codim I(o) \geq 2$
unless the pair ($C, \xi$) satisfies one of the following exceptional
conditions:
\par \noindent
(1) $o \in Sing X$ i.e. $\xi$ is not stable \par \noindent
(2) $C$ is hyperelliptic   \par \noindent
(3) there exists a double covering $\pi:C \rightarrow Y$ of an elliptic curve
and $\xi$ $=$ $\pi^*\eta
\otimes L$ where $L^2 = \omega_C \otimes det(\pi^* \eta^*)$ and $\eta$ is the
unique
(up to twisting by a degree zero line bundle) irreducible rank
two vector bundle of degree 1. \par \noindent
(4) $C$ is a smooth quartic curve in $\bold P^2$,  $\xi = T_{\bold
P^2}(-1)\otimes \Cal O_C(e)$,
($2e \sim 0$).\endproclaim \bigskip \noindent The next result is what we need
in the proof of the
main theorem: it is an immediate consequence of the previous statements 5.6 and
5.7.
\proclaim
{(5.8)THEOREM} Let $C$ be not hyperelliptic of genus $\geq 4$, $o$ a point of
$X-Sing X$, $o$
not the moduli point of the previous bundle $\pi ^* \eta \otimes L$ if $C$ is
bielliptic.
Assume $j$ is generically defined on $o \times T$ and moreover that $$ p_1
\cdot j : o \times T
\longrightarrow X $$ extends to a constant map, where $p_1: XÊ\times T
\rightarrow T$ is the first
projection.  Then $$ I_j \cap o \times T $$ has codimension $\geq 2$ in $T$.
\endproclaim

\noindent The first theorem we want to show is 5.7. Let $\xi$ be as above, we
need to consider one
closed set more: $$ H(\xi) = \lbrace l \in Pic^1(C) / h^0(\xi(l)) \geq 3
\rbrace \tag 5.9 $$
\noindent it is  clear that
$$
l \in H(\xi) \Longleftrightarrow -C_l \cup C_l \subset
\Theta_{\xi}
$$
\noindent where $C_l = \lbrace l-x, \quad x \in C \rbrace $
 \proclaim
{(5.10) LEMMA} $g-3 \leq dim H(\xi) \leq g-2 $ \endproclaim \demo {Proof} Let
$ \alpha : Pic^1(C) \times C \longrightarrow Pic^0(C)$ be the difference map,
consider
the natural projection \par \noindent  $ p: \alpha^* \Theta_{\xi}
\longrightarrow
Pic^1(C) $ and notice that $ p^*(l) = C_l \cdot \Theta_{\xi}$ for each $l \in
Pic^1(C)$. Then $p$ is generically finite on each irreducible component of
$\alpha^* \Theta_{\xi}$ and the locus of points $l$ where the fibre of $p$ is
not finite is $H(\xi)$:
therefore $dim H_{\xi}Ê\leq g-2$. On the other hand, if $ H(\xi) \cap
C^{(3)}-t$ is
not empty for any $t \in T$ it follows from transversality of general translate
and  the ampleness of
the cycle $C^{(3)}-t$ that $g-3 \leq dim H(\xi)$. It is known ([G]) that, for
at least one effective
divisor $d$ of degree 3, one has $h^0(\xi(t-d))Ê\geq 1$. Putting $l = d-t$ and
applying Riemann-Roch
it follows $h^1(\xi(-l)) = h^0(\xi(l)) \geq 3$. Hence $ d-t \in H(\xi) \cap
C^{(3)}-t \neq
\emptyset$.\enddemo \proclaim {(5.11) LEMMA} $V_{00}(\xi)$ and $V_m(\xi)$, $m =
0 \dots 3$, are
proper closed subsets of $T$, moreover: \par \noindent (1) $V_{00}(\xi)
\subset V_0(\xi)$, $codim V_0(\xi) \geq 2$,  $codim V_{00}(\xi) \geq 3$ \par
\noindent (2) $codim
V_1(\xi) \geq 1$ \par \noindent (3) $codim V_2(\xi) \geq 2$ if $\Theta_{\xi}$
is reduced
\par \noindent (4) $codim V_3(\xi) \geq 1$.\endproclaim
\demo {Proof} (1) observe that $t \in V_{00}(\xi)$ if and only if the surface
$ S_t = \lbrace
t-x-y,\quad x,y \in C \rbrace $ is in $\Theta_{\xi}$. Fix a point $e$ not in
$\Theta_{\xi}$ and
consider in $T$ the surface $S_e = \lbrace e+x+y, \quad x,y \in C \rbrace$.
Since $S_e$ is ample
and $S_e \cap V_{00}(\xi) = \emptyset$ it follows $codim V_{00}(\xi) \geq 3$.
Replacing the
surface $S_e$ by the curve $\Delta_e = \lbrace e+2x, \quad x \in C\rbrace$ and
repeating word by word
the same argument one shows $codim V_0(\xi) \geq 2$.\par \noindent (2): by
definition $t \in V_1(\xi)$
if and only if $ t = l+x $ where $x \in C$ and $l \in H(\xi)$. Considering the
sum map $ \beta: H(\xi)
\times C \longrightarrow T $, we have $V_1(\xi)= \beta(H(\xi) \times C)$. Since
$codim H(\xi)
\geq 2$ it follows $codim V_1(\xi) \geq 1$.  \par \noindent (3): by a theorem
of Laszlo [L3]
$h^0(\xi(t-2x)) \geq 2$ implies $t-2x \in Sing \Theta_{\xi}$. Therefore
$V_2(\xi) \subseteq
\beta (Sing \Theta_{\xi} \times \Delta)$ where $\beta$ is the sum map and
$\Delta \subset C(2)$  the
diagonal. The result follows immediately.   \par \noindent (4): observe that $t
\in
V_3(\xi)$ if and only if $t = 3x -l$ where $x \in C$ and $h^0(\xi(-l)) \geq 1$.
The latter condition
is equivalent to $h^0(\xi(l) \geq 3$ so that $l \in H(\xi)$. Consider the  map
$ \gamma : H(\xi)
\times C \longrightarrow T$ which is so defined: $\gamma(l,x) = 3x-l$; then
$V_3(\xi) = \gamma
((H(\xi) \times C)$. Since $codim H(\xi) \geq 2$ it follows $codim V_3(\xi)
\geq 1$. \enddemo
\proclaim {(5.12)LEMMA} Let $Z$ be an irreducible component of $V_{13}(\xi)$.
Assume $Z$ is not
contained in $W(\xi)$, then $$ codim Z \geq 2 $$ \endproclaim \demo {Proof} Let
$U = Z-(W(\xi) \cap
Z)$, we consider in $U \times C$ the closed set $$ \tilde U = \lbrace (t,x) \in
U \times C / \quad
h^0(\xi(t-x)) \geq 3 \quad \text {\it and \rm} \quad h^0(\xi(t-3x)) \geq 1
\rbrace $$ \noindent Under
the difference map  $ \alpha: \tilde U \longrightarrow Pic^1(C)$, ($\alpha(t,x)
= t-x$),  we have
$\alpha(\tilde U) \subseteq H(\xi)$ hence $\alpha(\tilde U)$ is at most
$g-2$-dimensional
by lemma 5.10. On the other hand the fibre of $\alpha $ is either the curve $C$
or a finite set. Let
$Y$ be any irreducible component of $\alpha (\tilde U)$ having maximal
dimension $g-2$: if, for $l$
general  in $Y$,  $\alpha^{-1}(l)$ is finite the result follows.  Let  $$ Y_n =
\lbrace l \in Y/
\quad h^0(\xi(l)) \geq n \rbrace $$\noindent it is obvious that $Y_3 = Y$,
moreover $Y_n$ is closed
in $Y$. We first assume that $Y_4$ is a proper subset so that $h^0(\xi(l)) = 3$
for $l$ general.
Consider the determinant map $$ d : \wedge^2 H^0(\xi(l)) \longrightarrow
H^0(\omega_C(2l)) $$ and
assume $\alpha^{-1}(l)$ is not finite, then $ \alpha^{-1}(l) $ $=$ $ \lbrace
(l+x,x), x \in C
\rbrace$ so that $$ h^0(\xi(l+x-3x)) = h^0(\xi(l-2x)) \geq 1 \forall x \in C$$
\noindent This is
impossible if $d$ is injective because then, in the 2-dimensional linear system
 defined by $Im d$ we
would have a pencil of divisors containing $2x$ for each $x$ in $C$. Hence
there exist  two
independent vectors $s_1,s_2 \in H^0(\xi(l))$ such that $d(s_1 \wedge s_2) = 0$
and they define  a
subline bundle $L$ of $\xi(l)$ with $h^0(L) \geq 2$, $deg L \leq g$. Since $L$
is also a  subline
bundle of $\xi(l+x)$ it follows $l+x \in W(\xi), \forall x \in C$: hence $l$
cannot be in $\alpha
(\tilde U)$, a contradiction.\par \noindent Secondly we assume $Y = Y_4$. Note
that, fixing a point
$x \in C$, $x+Y_5$ is contained in $V_{00}(\xi)$ which is at most
$g-3$-dimensional by the previous
lemma. Hence, generically on $Y$, $h^0(\xi(l)) = 4$. Assume $\alpha^{-1}(l)$ is
 not finite for a
general $l \in Y$ then  $$ h^0(\xi(l-2x)) \geq 1 \quad \text { for \it all \rm
$x \in C$ and moreover
$h^0(\xi(l)) = 4$}$$ \noindent to complete the proof we consider the surface $R
=\bold
P(\xi(l)^*)$ together with its tautological bundle $H$ (i.e. $p_*H = \xi(l)$)
and the rational map $$
f_H: R \longrightarrow \bold P^3 = \bold PH^0(\xi(l))^* $$ \noindent associated
to $H$, $R_x =
p^{-1}(x)$. Let $p:R \rightarrow C$ be the obvious projection, $R_x =
p^{-1}(x)$, we recall that
$h^0(\xi(l-mx)) = h^0(\Cal O_R(H-mR_x))$. There are two cases to be considered:
\par \noindent CASE
(1) $h^0(\Cal O_R(H-mR_x) \geq 2$ for all $x \in C$ and a fixed $m\geq 2$. \par
\noindent Observe that
the tangent map $(df_H)_z$ is not injective if $h^0(\Cal O_R(H-2R_x)) \geq 2$
and $z \in R_x$:
therefore $f_H(R)$ cannot be 2-dimensional in this case. Since $dim f_H(R) \leq
1$ there exists a
subline bundle $L$ of $\xi(l)$ with $h^0(L) = 4$ ($f_L(C) = f_H(R)$). Hence
$l+x \in W(\xi), \forall
x \in C$ and $l$ cannot  belong to $\alpha(\tilde U)$.\par \noindent CASE (2)
$dim f_H(R) = 2$ and
$h^0(\Cal O_R(H-2R_x)) = 1$ for a general $x \in C$. \par \noindent Assuming
this and keeping $x$
general  we observe that there exists a unique element $H_x \in \mid H-2R_x
\mid $, moreover $ \mid
H_x - R_x \mid $ is (generically) empty: otherwise, fixing the point $x$, we
would have $3x-l \in
V_{00}(\xi)$ for $l$ general, hence $dim Y_4 \geq g-3$. From these two remarks
it follows that the
pencil $\mid H - R_xÊ\mid $ has a unique base point $b(x)$  on the fibre $R_x$.
This defines the
holomorphic section $$ b: C \longrightarrow R $$ \noindent sending $x$ to
$b(x)$ and the curve $B =
b(C) $. Let $i(x)$Êbe the intersection index  of $B$ and $H$ at $b(x)$: assume
$\mid H \mid$ has no
base points on $R_x$ then, by the construction of $B$,  $$ i(x) \geq 2 \quad
\text { if and only if
$H$ contains $R_x$}Ê$$ \noindent This implies that either $f_H(R_x)$ is a
tangent line to the curve
$f_H(B)$ or $f_H(R)$ is a cone of vertex the point $f_H(B)$. In any case fix a
general pencil $P
\subset \mid H \mid$ and consider its restriction $P_B$ to $B$: if $d$ is the
degree of the moving
part of $P_B$ one computes from the previous equivalence and Hurwitz formula:
$2d -  H^2 \leq 2- 2g
$, hence $d \leq 1$ because $H^2 = deg \xi(l) = 2g$. Since this is impossible
$P_B$ has no moving
part and an element of $P$ must contains $B$; moreover $dim \mid H-B \mid \geq
3$ and $f_H(R)$ is a
cone. Therefore there exists a subline bundle $L$ of $\xi(l)$ with $h^0(L) \geq
3$ and one completes
the proof as in case (1).\enddemo \noindent The next lemma is standard, we omit
for brevity its proof
\proclaim {(5.13) LEMMA} (1) Let $ Y \subset Pic^n(C) $ be a  closed subset,
$\alpha:
 C^{(n)} \times Y \longrightarrow J = Pic^0(C) $
the difference map, $Z = \alpha (C^{(n)} \times Y)$. If $dim Z <
g$ then $\alpha$ is generically finite onto $Z$. In particular: \par \noindent
(2) Let $ Y
=\lbrace l \in Pic^n(C) / C^{(n)}-l \subseteq Z,  \rbrace $ where $Z$ is a
proper closed subset of
$Pic^0(C)$. Then $ dim Y \leq dim(Z)-n $. \endproclaim
\noindent The next result is an application of Martens' theorem: \proclaim
{(5.14) THEOREM} For a semistable rank two vector bundle $\xi$ of determinant
$\omega_C$Êlet
$$\text {(5.15) $ W^1_n(\xi,i) = \lbrace l \in Pic^i(C) / $ there exists $L \in
Pic^n(C)$ such that
$h^0(\xi(l) \otimes L^{-1}) \geq 1$, $ h^0(L) \geq 2 \rbrace$}    $$ assume $0
\leq i \leq n < g$
then: \par \noindent (1) if $C$ is not hyperelliptic  $ dim W^1_n(\xi,i) \leq
g+i-4 $ \par \noindent
(2) if $C$ is hyperelliptic
 $ dim W^1_n(\xi,i) = g+i-3 $ for $i \leq 3$ \par \noindent
 Assume $0 \leq i \leq n  = g$ then: \par \noindent
 (3) $dim W^1_g(\xi,i) \leq g+i-3$
\endproclaim
\demo {Proof}ÊLet $W^1_n \subset Pic^n(C)$ be the Brill-Noether locus of line
bundles
$L$ with $h^0(L) \geq 2$. We can construct in $ W^1_n \times Pic^i(C)$ the
closed
subset
$$
B = \lbrace (L,l) \in W^1_n \times Pic^i(C) / h^0(\xi(l)\otimes L^{-1}) \geq
1Ê\rbrace
$$
considering the two  projections
$$
\CD
{W^1_n} @<{p_1}<< B @>{p_2}>> {Pic^i(C)} \\
\endCD
$$
it is clear that
$$
W^1_n(\xi,i) = p_2(B)
$$
On the other hand, fix
$
L \in p_1(B)
$
and consider the fibre
$$
F_L = p^{-1}_1(L)
$$
Let $l \in F_L$, $m = L(-l) $ then $h^0(\xi(-m)) \geq
1$ so that $h^0(\xi(d-m))\geq 1$ for all $d \in C^{(n-i)}$ and
$$ C^{(n-i)}-m \subseteq \Theta_{\xi}$$
In particular $F_L \subseteq \lbrace m \in Pic^{(n-i)} / C^{n-i}-m
\subseteq \Theta_{\xi}Ê\rbrace$, hence $ dim(F_L) \leq (g-1) -
(n-i) $ by the previous corollary; on the other hand, by Martens theorem,  $dim
W^1_n \leq n-3$ if
$C$ is not hyperelliptic. Therefore $$ dim W^1_n(\xi,i) \leq dim B \leq dim
W^1_n + dim F_L \leq
g+i-4   $$ We leave as an exercise the hyperelliptic case. The proof of the
case
$n = g$  is exactly the same: of course in this case $dim W^1_g = g-2$ so that
$dim W^1_g(\xi,i) \leq g+i-3$ . \enddemo
Finally, from a theorem of Lange and Narasimhan [LN], we obtain \proclaim
{(5.16)
PROPOSITION} $Codim W(\xi) \geq 2$ unless: \par \noindent (1) $\xi$ is not
stable
\par \noindent (2) $C$ is hyperelliptic   \par \noindent
(3)there exists a double covering $\pi:C \rightarrow Y$ of an elliptic curve
and $\xi$ $=$ $\pi^*\eta
\otimes L$ where $L^2 = \omega_C \otimes det \pi^* \eta$ and $\eta$ is the
unique irreducible rank
two vector bundle of degree 1, (up to twisting by a degree zero line bundle).
\par \noindent
(3) $C$ is a smooth quartic curve in $\bold
P^2$,  $\xi = T_{\bold P^2}(-1)\otimes \Cal O_C(e)$, $2e \sim 0$. \endproclaim
\demo {Proof} Using the previous notations we have $W(\xi) = \bigcup
W^1_n(\xi,2), \quad 2 \leq n \leq g $ \par \noindent
By the previous theorem 5.14 $codim W^1_n(\xi,2) \geq 2$ if $C$ is not
hyperelliptic
and $n \leq g-1$. So we have only to consider the case $n = g$: by definition
$t \in W^1_g(\xi,2)$ iff there exists $L_t \in W^1_g$ such that $L_t(-t)$ is a
subbundle of $\xi$.
A simple dimension count shows that if $codim W^1_g(\xi,2) = 1$ then $\xi$ has
(at least) a
1 dimensional family  of such subline bundles of degree $g-2$.
If $\xi$ is stable
 they are of maximal degree$g-2$, hence by a theorem of
Lange-Narasimhan $\xi$ ([LN] thm. 5.1) , $\xi$ is as in cases (3) or (4). If
$\xi$
is semistable not stable it is easy to see that $codim W(\xi) = 1$.
\enddemo  \bigskip \noindent \bf PROOF OF THEOREM 5.7\rm:Êthis is now an
immediate consequence of the previous results: let $o = [\xi]Ê\in X $, by lemma
5.12 and 5.13 the only
component of  $I(o)$ having possibly codimension 1 is the "Brill-Noether locus"
$W(\xi)$. On the
other hand, by the previous theorem and proposition 5.16, one has, with the
prescribed exceptions,
$codim W(\xi) = 2$. Hence the same holds for $I(o)$. \bigskip \noindent In
order to show  theorem
5.6 we need some preparation. In particular we want to give the following
\bigskip \noindent \bf
ALTERNATIVE DEFINITION OF $j$Ê\rm: \bigskip \noindent Let  $(o,t) \in X \times
T$,  $$ o = [\xi] $$
we  denote by $$ e: H^0(\xi(t)) \otimes \Cal O_C \longrightarrow \xi(t) $$ the
evaluation map and by
$$ \Cal U^* \tag 5.17$$ \noindent the image sheaf $Im(e)$. Then we consider the
two projections  $$
\CD C @<{p_1}<< { C \times C} @>{p_2}>> C \\ \endCD \tag 5.18$$ and the
diagonal  $$ \Delta \subset C
\times C  $$ \noindent Tensoring
$$
0 \longrightarrow \Cal O_{C\times C}(-2\Delta) \longrightarrow \Cal O_{C\times
C}(-\Delta)
\longrightarrow \Cal O_{\Delta}(-\Delta) \longrightarrow 0
$$
by
$$
\Cal F= p_1^* \Cal U^* \tag 5.19
$$
and applying the $p_{2*}$ functor, we obtain the long exact sequence
\par \noindent

$$
\CD
0 @>>>Ê{p_{2*} \Cal F(-2\Delta) } @>>> {p_{2*} \Cal F(-\Delta)} @>>> \\
@>>> {p_{2*} \Cal F \otimes \Cal O_{\Delta}(-\Delta) } @>>> {R^1p_{2*} \Cal
F(-2\Delta)}Ê@>>> { R^1p_{2*} \Cal F(-\Delta)} @>>> 0 \\
\endCD \tag 5.20
$$Ê
\noindent Before of using it we fix the following open set in $X \times T$
\proclaim {(5.21)DEFINITION}
$$
U_{j'} = \lbrace (o,t) \in X \times T / \quad  \text { t  is not in  $V_0(\xi)
\cup
V_{13}(\xi)$,} \quad o = [\xi] \rbrace $$
\endproclaim
\noindent
let $ (o,t) \in U_{j'}$, since $t$ is not in $V_{00}(\xi)$ it follows:
\bigskip \noindent
(1) $h^0(\xi(t)) = 4$ \par \noindent
(2) the evaluation map $e: V \otimes \Cal O_C \rightarrow \xi(t)$ is
generically surjective, where
$V$ $=$ $H^0(\xi(t))$. \bigskip \noindent This defines the exact sequence
$$\CD
0 @>>> \Cal U @>{e^*}>> {V^*\otimes \Cal O_C} \longrightarrow \overline {\Cal
U} @>>> 0 \\
\endCD \tag 5.22
$$
\noindent
where both
$$
\Cal U = Im(e^*) = Im(e)^* \quad \text {and} \quad \overline {\Cal U}
$$
are rank two vector bundles. Considering the Grassmannian $G$  of 2-spaces in
$V^*$
and the Gauss map 1.5 $$  g: C \longrightarrow G $$ of the pair $(V,\xi(t))$ it
turns out that the
previous sequence is the pull back by $g$ of the universal/quotient bundle
sequence on $G$ (see
section 1) $$
0 \longrightarrow \Cal U_V \longrightarrow V^* \otimes \Cal O_C \longrightarrow
\overline {\Cal U}_V
\longrightarrow 0 $$

\proclaim {(5.23)LEMMA} If $(o,t) \in U_{j'}$ the sequence (5.20) becomes
$$
\CD
0 Ê@>>> {\overline {\Cal U}^*} @>u>> {\Cal U^* \otimes \omega_C} @>>>
{R^1p_{2*} \Cal
F(-2\Delta)}Ê@>>> { R^1p_{2*} \Cal F(-\Delta)} @>>> 0 \\
\endCD
$$Ê
In particular this induces the exact sequence
$$
\CD
0 Ê@>>> {\overline {\Cal U}^*} @>u>> {\Cal U^* \otimes \omega_C} @>>> {\Cal D}
@>>> 0 \\
\endCD \tag 5.24
$$
\noindent where $\Cal D = Coker u $ is an $O_C$-module of finite length. $\Cal
D$ is supported
exactly on points $x$ such that $h^0(U^*(-2x)) \geq 1$.\endproclaim
\demo {Proof} Let $z \in C$: by definition, tensoring (5.20) by $\Cal
O_{C,z}/m_z$, we
obtain the long exact sequence
$$\CD
0 @>>>{H^0(\Cal U^*(-2z))} @>>> {H^0(\Cal U^*(-z))} @>u_z>> {\Cal U^*_x\otimes
\omega_C} @>>> \\
@>>> {H^1(\Cal U^*(-2z))}Ê@>>> {H^1(\Cal U^*(-z)} @>>> 0 \\
\endCD
$$
\noindent Since $\Cal U^*$ is globally generated and $H^0(\Cal U^*) = V$ it
follows
that $H^0(\Cal U^*(-z))$ is of constant dimension 2. Since $H^0(\Cal U^*(-z))$
is the fibre of
$p_{2*} \Cal F(-\Delta)$ at $z$ this latter sheaf is a rank two vector bundle.
Furthermore the
universal/quotient bundle sequence $$ 0 \rightarrow \overline {\Cal U}^*
\rightarrow V \otimes \Cal
O_C \rightarrow \Cal U^* \rightarrow $$ yelds the canonical isomorphims
$\sigma_z: \overline {\Cal
U}_z \rightarrow H^0(\Cal U^*(-z))$, $z \in C$. This implies that $p_{2*} \Cal
F(-\Delta)$ is
isomorphic to $\overline {\Cal U}^*$.\par \noindent Note that $\Cal
O_{\Delta}(-\Delta) \cong
\omega_C \cong p_2^* \omega_C$. Therefore, by projection formula,  $$ p_{2*}
\Cal F \otimes \Cal
O_{\Delta}(-\Delta) = p_{2*}(p^*_1\Cal U^*\otimes \Cal O_{\Delta}) \otimes
p^*_2\omega_C) = \Cal U^*
\otimes \omega_C $$\noindent Since $t$ is not in $V_0(\xi)$ we have $h^0(\Cal
U^*(-2z) =
h^0(\xi(t-2z)) = 0$ for a general $z$. This implies that $p_{2*} \Cal
F(-2\Delta)) = 0$ and that
$Coker u$ is supported on points $z$ such that $h^0(\Cal U^*(-2z)) > 0$.
\enddemo Since $\Cal U^* =
Im(e)$  we have the exact sequence  $$\CD 0 @>>> {\Cal U^*} @>>> {\xi(t)} @>>>
{\Cal C} @>>> 0 \\
\endCD \tag 5.25 $$  \bigskip \noindent with
$$
\Cal C = Coker (e) \tag 5.26
$$
\par \noindent Both
$$
D = Supp \Cal D \quad \text {and } \quad Z = Supp \Cal C
$$ are divisors in $C$, in particular
$$
\Cal C = \Cal O_Z \quad , \quad \Cal D = \Cal O_D
$$
\noindent We want to point out that
$$
Supp Z \cap Supp D = \emptyset \quad \text {so that} \quad \Cal D \otimes \Cal
C = 0 \tag 5.27
$$
this is due to our choice of $(o,t)$ in $U_{j'}$ because then $t$ is not in
$V_{13}(\xi)$. Another
remark is that, whenever $\xi(t)$ is globally generated and $\overline {\Cal
U}$ semistable, then
$det \overline {\Cal U}(-t)$ $\cong \omega_C$ and moreover
$$
j(o,t) = (\overline o,t)
$$
for the moduli point $\overline o$ of $\overline {\Cal U}(-t)$: this just
follows from the definition
of $j$. We want to use the previous constructions to give a partial extension
of $j$ to points
$(o,t)$ such that $\xi(t)$ is not globally generated.  With this purpose we
consider the exact
sequence $$\CD
0 @>>> {\xi(t)}^* @>>> {\Cal U  } @>{\epsilon}>> {\Cal C} @>>> 0 \\ \endCD $$
\noindent which is the
dual of 5.25 and note that the map $\epsilon$ is an element of $$ \Cal U^*
\otimes \Cal C = Hom(\Cal
U, \Cal C)$$ \noindent Then we tensor by $\Cal C$ the sequence
$$
\CD
0 Ê@>>> {\overline {\Cal U}^*\otimes \omega_C^*} @>u>> {\Cal U^* } @>>> {\Cal D
\otimes {\omega_C}^*}
@>>> 0 \\
\endCD
$$
\noindent constructed in 5.24. Since $\Cal C \otimes \Cal D  = 0$  we obtain
$$ \overline {\Cal
U}^*Ê\otimes \Cal C \otimes \omega_C^* \longrightarrow \Cal U^* \otimes \Cal C
\longrightarrow \Cal
C \otimes \Cal D \otimes \omega_C^* = 0 $$ \noindent and hence a natural \it
isomorphism \rm  $$
\sigma: \overline {\Cal U}^*Ê\otimes \Cal C \otimes \omega_C^* \longrightarrow
\Cal U^* \otimes \Cal
C  $$ Let
$$
\overline {\epsilon} = \sigma^{-1}(\epsilon)
$$
then $\overline {\epsilon} \in Hom(\overline {\Cal U} \otimes \omega_C, \Cal
C)$. Since
$\epsilon$ is an epimorphism it is easy to see that $\overline {\epsilon}$ is
an epimorphism too.
This  defines the rank two vector bundle $$ Ker \overline {\epsilon} $$
\noindent such that $det (Ker \overlineÊ{\epsilon}) = det \xi(t) \otimes
\omega_C^{2}$. Then we
fix the following \proclaim {(5.28)DEFINITION} $Ker \overline {\epsilon}
\otimes \omega_C^*$ will be
always denoted by $ \overline {\phi}_t $ \endproclaim \noindent  and briefly
sketch the conclusion:
the construction of $Ker \overline {\epsilon}$ easily extends from one pair to
a family
of pairs $(\xi,t)$ such that $([\xi],t) \in U_{j'}$. Then it is standard to
deduce that there
exists a rational map $$ j': X \times T \longrightarrow X \times T \tag 5.29 $$
which is so defined on
$U_{j'}$: $$
\forall ([\xi],t) \in U_{j'} \quad  j'([\xi],t) = ([\overline {\phi}_t],t)
$$
provided $\overline {\phi}_t$ is semistable. From the definitions of $j$ and
$j'$ it follows
$$
j(o,t) = j'(o,t) \quad \text {on $U_{j'} \cap U_j$}
$$
where $U_j$ is the open set defined in section 2. Therefore $j'$ and $j$ define
the same birational
involution.  \bigskip \noindent Let us recall two elementary facts: \bigskip
\noindent
(5.29) REMARKS\par \noindent
(1) If $\xi(t)$ is globally generated $\Cal C = Coker(e) = 0$. Hence $\epsilon
= \overline
{\epsilon} = 0$ and $\overline {\phi}_t(t) = \overline {\Cal U}$ $=$ $Ker
(e)^*$.\par \noindent
(2)We have obvious canonical inclusions $$ V^* \subset H^0(\overline {\Cal U})
\subseteq H^0(\overline
{\phi}_t(t))^* $$
\noindent using the standard arguments of section 1 and the exact sequence
$$
0 \longrightarrow \Cal U \longrightarrow V^* \otimes \Cal O_C \longrightarrow
\overline {\Cal U}
\longrightarrow 0
$$
\noindent it is easy to see that the pairs $(V,\xi(t))$ and $(V^*,\overline
{\phi}_t(t))$
define the same quadric $Q$ containing the embedded curve $C_t \subset \bold
P^{g+2}_t$:
$$
Q = Q(V,\xi(t)) = Q(V^*, \overline {\phi}_t(t))
$$

\bigskip
\noindent \proclaim {PROOF OF THEOREM 5.6} \endproclaim
\demo {Proof} Let $o$ $=$ $[\xi]$, at first we assume  that $\xi$ is \it stable
\rm. The assumption
of the theorem we want to prove is that
$$
j(o,t) = (\overline o,t)
$$
for a general $t$ and a fixed point
$$
\overline o = [\overline {\xi}]
$$
\noindent As it is easy to show (see next lemma), for each point $u$
corresponding to a
semistable not stable bundle and $t$ general $j(u,t) = (u,t)$; therefore
$\overline {\xi}$ is stable too. Let $I(o,\overline o) = I(o) \cup I(\overline
o)$ be the special
set of $\lbrace o, \overline o \rbrace$ defined in 5.5, we must show that the
open
subset
$$
U= T-I(o,\overline o)
$$
is in the domain of $j$. Let $t \in U$, from the previous construction of $j'$
we already know that
$$
j'(o,t) = j(o,t) = ([\overline {\phi}_t],t)
$$
as soon as the rank two vector bundle $\overline {\phi}_t$ is semistable. To
show that
$\overline {\phi}_t$ is semistable we fix  a Poincar\'e bundle on $C \times T$
and construct in a
standard way a sheaf
$$
\Cal R \longrightarrow C \times T
$$
\noindent such that
$$
\Cal R_t = \Cal R \otimes \Cal O_{C \times t} = \overline {\phi}_t(t)
$$
\noindent for all $t \in U$. On the other hand, for a \it general \rm $t \in
U$, we have
$$
\overline{\xi}(t) \cong \overline {\phi}_t(t)
$$
\noindent To see this take $t$ in the complement of $I(o,\overline o) \cup
V_1(\xi)
\cup V_1(\overline {\xi})$ and  such that $(o,t)$ is in the domain of$j$. Under
this choice both
$\xi(t)$ and $\overline {\xi}(t)$ are globally generated and moreover the
sequence
$$
0 \longrightarrow \overline {\xi}(t)^* \longrightarrow H^0(\xi(t)) \otimes \Cal
O_C \longrightarrow
\xi(t) \longrightarrow 0
$$
\noindent is exact. Since $\xi(t)$ is globally generated it is clear that the
construction of
$\overline {\phi}_t(t)$ gives exactly $\overline {\xi}(t)$, (cfr. remark
5.29(1)). By
semicontinuity the condition $ \overline \xi(t) \cong \Cal R_t $ for a general
$t$ implies $$
h^0(\Cal R_t \otimes \overline {\xi}(t)^*) \geq 1 $$   for all $t \in T$.
Equivalently, there exists
a  non zero morphism  $$ \rho_t : \Cal R_t \longrightarrow \overline {\xi}(t)
$$ \noindent for each
$t \in T$. Assume $ t \in U $ and  $\Cal R_t = \overlineÊ{\phi}_t(t) $
unstable, then $\rho_t$ cannot
be an isomorphism; since $\overlineÊ{\phi}_t$, $\overline {\xi}$ have the same
slope and the latter
is stable we have an exact sequence $$ 0 \longrightarrow \Cal A \longrightarrow
 \overline
{\phi}_t(t) \longrightarrow \Cal B \longrightarrow 0 $$ where $$ \Cal A = Ker
(\rho_t) \text { and }
\Cal B = Im (\rho_t) $$ are line bundles and $\Cal B$ is a subline bundle of
$\overline {\xi}(t)$ of
degree $\leq g$. To obtain a contradiction we  show that  $$ h^0(\Cal B) \geq 2
$$ \noindent because
then $h^0(\Cal B) \geq 2$ $\Rightarrow$ $t \in W(\overline {\xi}) \subset
I(o,\overline o)$
$\Rightarrow$ $t$ is not in $U$.  By construction (see section 1 and remark
5.29) we have   $$
q(V,\xi(t)) = q = cq(V^*,\overlineÊ{\phi}_t(t))
$$
\noindent ($c \in \bold C^*$) for the pairs $(V,\xi(t))$ and $(V^*,\overline
{\phi}_t(t))$; the
quadratic form $q$ is vanishing on the embedded curve $C_t$. If $h^0(\Cal B)
\leq 1$ then $H^0(\Cal
A)$ cuts on $V^*$ a subspace of codimension $\leq 1$ and hence $q$ is
identically zero by prop. 1.11.
Since $q$ is also defined by $(V,\xi(t))$ the same proposition implies that
$\xi(t)$
has a subline bundle $\Cal C$ with $h^0(\Cal C) \geq 3$: this is impossible
because $t$ is in $U$,
hence not in $W(\xi)$. Therefore $h^0(\Cal B) \geq 2$. It remains to consider
the semistable
case: this is done in the next lemma. \enddemo  \proclaim {(5.30) LEMMA} Let $o
\in Sing(X)$ then:
\par \noindent (1) $I_j \cap o \times T \subseteq I(o)$ \par \noindent (2)
$j(o,t) = (o,t)$ for $t
\in T-I(o)$ \endproclaim \demo {Proof} Since $o \in Sing X$ $o$ is the moduli
point of a family of
(S-equivalent, see [B1]) semistable not stable vector bundles. Exactly one of
them is split: let $$
\xi = M \oplus N $$ such a bundle: if $t \in V_1(\xi)$ then $\mid M(t) \mid$
(or $\mid N(t) \mid $)
has a base point $x$ and $L = M(t-x)$ is a subline bundle of $\xi(t)$ with
$h^0(L) \geq 2$, $deg L
\leq g$. Hence $t \in W(\xi)$. Conversely, if $t \in W(\xi)$, it is easy to see
that $\xi(t)$ is not
globally generated so that $t \in V_1(\xi)$. Therefore $V_1(\xi)$ $=$ $W(\xi)$
so that $$I(o) = W(\xi)
\cup V_0(\xi)$$ \noindent If $t \in T-I(o)$, one has immediately the exact
sequence $$ 0 \rightarrow
M(t)^*\oplus N(t)^* \rightarrow H^0(M(t))\otimes \Cal O_C \oplus H^0(N(t))
\otimes \Cal O_C
\rightarrow M(t) \oplus N(t) \rightarrow $$ which implies $j(o,t) = (o,t)$.
This completes the proof.
\enddemo
\noindent (5.31) REMARK Using the previous lemma and our main theorem one could
deduce that
$I_j \cap o \times T = W(\xi)$ if $ [\xi] \in Sing X$ and $C$ is not
hyperelliptic. Note
also that $W(\xi)$ is a divisor if $[\xi] \in Sing X$.
\proclaim {\bf 6. An application of Rigidity
lemma \rm} \endproclaim \bigskip \noindent In this section we fix $$ X=  \text
{integral quasi
projective variety} \tag 6.1 $$ and  $$\Cal L= \text {ample, globally generated
line bundle on $X$}
\tag 6.2 $$
 so that the  map associated to $\Cal L $
$$
\theta: X \longrightarrow \bold P^N=\bold P(H^0\Cal L)^*  \tag 6.3 $$
is a finite morphism. We want to check whether  $\Cal L$ is very ample by
composing $\theta $  with the elements of  a \it complete \rm family of linear
projections of $\bold P^N$. More precisely this means that we  fix  $$ T= \text
{integral projective  variety} $$
and a non zero map of vector bundles over $T$
$$
\lambda : H^0(\Cal L)^*\otimes \Cal O_T \longrightarrow \Cal Q  \tag 6.4
$$
Denoting  by $$\overline {\lambda} : \bold P^N \times T \longrightarrow \bold
P(\Cal Q)
\tag 6.5 $$ the induced map of projective bundles  and composing it with $$
\theta \times id_T : X \times T \longrightarrow \bold P^N \times T \tag 6.6 $$
we obtain a  rational map
$$ \Phi = \overline {\lambda} \cdot (\theta \times id_T) :  X \times T
\longrightarrow  \bold P(\Cal Q)\tag 6.7  $$ Clearly, for each $t \in T $, the
restriction  $$\Phi_t : X \times t \longrightarrow \bold P(Q_t)
$$
of $\Phi$ to $X \times t$ is just the composition of $\theta $ with the linear
projection
$$\overline {\lambda}_t : \bold P^N \rightarrow \bold P(\Cal Q_t).$$
In principle, when $X$ is projective, one can try to show  the very ampleness
of $\Cal L$ by
checking if any  pair of points (of tangent vectors)  can be separated by at
least  one
$\Phi_t$. Under some suitable assumptions such a family of rational maps
$\Phi_t$
can also be used to show that
$$ deg(\theta)=1 \Longrightarrow \text {$\Cal L$ very ample.}$$ In this section
we show
this result under some assumptions which are "ad hoc" constructed for the
case $X $ = $SU_2(C)$, $\Cal L $ = generalized theta divisor.  Essentially,
every
result here is a consequences of the well known \proclaim {RIGIDITY LEMMA}  Let
$ \Phi: X
\times T \longrightarrow Z $ be a morphism. Assume that there exists $o \in X$
such that $$ \Phi/ (o  \times T) \quad \text {is constant} $$ then $\Phi$
factors
through the first projection $p_1: X \times T \rightarrow T$
 i.e. $$ \Phi=f \cdot p_1 $$ where $f: X \rightarrow Z$ is a morphism.
\endproclaim \noindent As a slight extension of the lemma we have
\proclaim {(6.8)PROPOSITION } Let $\Phi: X \times T \rightarrow Z$ be only a
rational
map, $I_{\Phi}$ its indeterminacy locus. Assume there exists $o \in X$ such
that \par
\noindent  (1) $codim_T(I_{\Phi} \cap  o \times T) \geq 2$  \par \noindent (2)
$\Phi$
restricted to $  o  \times T $ extends to a constant map. \par \noindent Then
$\Phi$
factors through the first projection $p_1:X \times T \rightarrow X$ i.e.:
$$
\Phi= f \cdot p_1
$$
where $f: X \rightarrow Z$ is a rational map.\endproclaim \demo
{Proof} Let  $$A = \lbrace x \in X | codim_T(I_{\Phi} \cap  x \times T) \geq 2
\rbrace$$ by assumption (1) and dimension theory $A$ is a non empty open set.
Let us
show that $\Phi$ restricted to $ x \times T$ extends to a constant map for all
$x \in
A$: take any two points $$ t_1,t_2 \in T $$ such that $(x,t_1),(x,t_2)$ are not
in $I_{\Phi}$.
Since $T$ is projective and $codim_T (I_{\Phi} \cap   x  \times T)$  $\geq 2$,
$codim_T(I_{\Phi} \cap o \times T)$ $\geq 2$  there exists a \it projective \rm
curve
$\Gamma \subset T$ containing $t_1,t_2$ and such that
 $I_{\Phi} \cap (x \times \Gamma) = I_{\Phi} \cap (o \times \Gamma) = \emptyset
$. Now
consider the (non empty) open set $A_{\Gamma} = \lbrace x \in X / I_{\Phi} \cap
x \times
\Gamma = \emptyset \rbrace $:  by assumption (2) and the Rigidity Lemma applied
to $\Phi$
restricted to $A _{\Gamma} \times B$, we  obtain $$ \Phi(x,t_1)=\Phi(x,t_2)$$
Since $t_1,t_2$ are
choosen on an open set of $T$ it follows
 that $\Phi/ x \times T$ is a constant map. Then we can consider the function
 $$ f: A \longrightarrow Z $$
sending $x \in A$ to the point $\Phi(x \times T)$: clearly $f$ is a rational
map from
$X$ to $ Z$ and $\Phi = f \cdot p_1$. \enddemo Let us point out a special
corollary of this: \proclaim {(6.9)PROPOSITION  } Let $j: X \times T
\rightarrow X
\times T$ be a birational involution, $I_j$ its indeterminacy locus, $p_1: X
\times T
\rightarrow X$ the natural projection. Assume there exists $o \in X$ such that:
\par \noindent (1) $ codim_T(I_j \cap o \times T) \geq 2 $ \par \noindent (2)
$p_1
\cdot j$ restricted to $o \times T$ extends to a constant map \par \noindent
then
there exists a birational involution $f: X \rightarrow X$ which makes
commutative the
diagram $$ \CD {X \times T} @>>j> {X\times T} \\ @VV{p_1}V @VV{p_1}V \\ X @>>f>
X \\
\endCD $$ \endproclaim \demo {Proof} Consider the rational map $\Phi = p_1
\cdot j : X
\times T \rightarrow X$ and observe that $I_j$ contains the indeterminacy locus
of
$\Phi$; then apply to $\Phi$ proposition 6.8:   by assumptions (1) and (2)
$\Phi$ factors
through $p_1$ and a rational map $f: X \rightarrow X$; therefore $f \cdot p_1$
= $p_1
\cdot j$: this gives the required commutative diagram. In particular $f \cdot f
\cdot
p_1 $ = $ f \cdot p_1 \cdot j$ = $p_1 \cdot j \cdot j = p_1 $ so that $f \cdot
f =
id_X$ and $f$ is a birational involution. \enddemo
Now we apply the previous propositions  to the special situation which is
relevant for this paper; therefore let $\Phi: X \times T \rightarrow \bold
P(\Cal Q)$ be
the rational map constructed in 6.7 we make the following  \proclaim
{(6.10)ASSUMPTION } $\Phi$ has degree two onto its image. \endproclaim

\noindent Moreover, given such a $\Phi$, we will always use the following
\proclaim {(6.11) NOTATIONS} \par \noindent
(1) $j: X \times T \rightarrow X \times T$ =: the involution induced by $\Phi$
\par
\noindent (2) $I_j$ =: the indeterminacy locus of $j$ \par \noindent
(3) $B_j$ =: the maximal $j$-invariant open set such that $j/B_j$ is a
biregular
isomorphism. \endproclaim \noindent   As usual we say that $\theta$ is not an
embedding at $x \in X$ if $\theta$ is not injective at $x$ or if the tangent
map $d\theta_{x} : T_{X,x} \rightarrow T_{\theta(X),\theta(x)}$  is not an
isomorphism. Let $$ N(\Cal L) = \lbrace x \in X / \text {$\theta$ is not an
embedding at $x$ \rbrace } \tag 6.12 $$ then \proclaim {(6.13)THEOREM } Let
$\theta: X \longrightarrow \bold P^n$ be the morphism associated to $\Cal L$.
Assume: \par \noindent (1) $\forall x \in N(\Cal L) \quad codim_T(I_j \cap x
\times T) \geq 2 $ \par \noindent   (2) $\forall x \in N(\Cal L)$, $p_1 \cdot
j$ restricted to $x \times T$ extends to a constant map \par \noindent Then,
if $\Cal L $ is not very ample, there exists a non identical birational
involution $ f: X \longrightarrow X $ such that $$ j = f \times id_T$$
Furthemore, assume in addition that:\par \noindent (3) $X$ is normal, $Pic (X)
\cong \bold Z$ and the restriction $Pic (X) \rightarrow Pic (X-Sing X) $ is an
isomorphism \par \noindent (4) $codim_XY \geq 2$, where $Y = X-p_1(B_j)$ \par
\noindent Then, if $\Cal L$ is not very ample, there exists a projective
involution $\overline f$ of $\bold P^n$ such that   $$ \theta \cdot f =
\overline f \cdot \theta $$
 \endproclaim
\demo {Proof} If $\Cal L$ is not very ample
 $N(\Cal L)$ is not empty. Then, by proposition 6.9, there exists a
birational
involution $f$ on $X$ such that $p_1 \cdot j = f \cdot p_1$ where $p_1$ is the
first projection of $X \times T$. Since $j(x,t)$ is the identity in $t$ it
follows $ j = f \times id_T$ and in particular $f$ is non identical.\par
\noindent
Since $f$ = $j \times id_T$ the open set $X-Y$ is $f$-invariant and $f$ is
biregular on it. Let $Z = Y \cup Sing (X)$: from assumption (4) and the
normality
of $X$ it follows that $Z$ has codimension $\geq 2$ in $X$. Consider the open
set
$ U = X - Z $ from assumption (3) it follows that the
restriction  $Pic(X)  \rightarrow Pic(U)$ is an isomorphism. Clearly $U$ is
$f$-invariant and $f/U$ is a \it biregular \rm involution of $U$. Hence, in
particular, $(f/U)^* \Cal L/U$ $ \cong \Cal L/U$.  Since $X-U$ has codimension
$\geq 2$ the restriction $H^0(\Cal L) \rightarrow H^0(\Cal L/U)$ is an
isomorphism too. This implies that $f/U$  induces on $\bold P^N$
=Ê$\bold PH^0(\Cal L)^*$ a projective involution $\overline f$ such that $
\theta \cdot f = \overline f \cdot \theta$. \enddemo \par \noindent (6.14)
REMARK The previous theorem is only a  corollary to proposition 6.9,
nevertheless it plays an essential role for showing the very ampleness of the
generalized theta divisor. Let us see why:  of course if we  produce a subset
$S
\subset X$ such that $\theta(S)$ spans $\bold P^N$ and $f/S$ = $id_S$, then
$\overline f$ is the identity and  $deg (\theta) \geq 2$. In the case $X$ =
$SU_C(2)$, $\Cal L$ = generalized theta divisor we will see that such an $S$ is
$Sing X$;  therefore  we will be able to conclude $\Cal L$ not very ample
$\Longrightarrow$ $deg(\theta) \geq 2$ $\Longrightarrow$ $C$ hyperelliptic,
where the second implication is  Beauville's result (0.4) (iii).
\proclaim {7. Proof of the main theorem }\endproclaim
 \noindent
 Let $\theta : X \longrightarrow \mid 2\Theta \mid $ be the morphism associated
to the generalized
theta divisor, we want to show our main theorem (0.3) that is: \bigskip
\noindent
\it if $C$ is not hyperelliptic \par \noindent
(1) $\theta$ is injective \par \noindent
(2) $d\theta_x$ is injective $\forall x \in X-Sing X$. \rm \bigskip \noindent
Since everything is known for genus $g \leq 3$ we assume $C$ not hyperelliptic
of genus $g \geq 4$,
then we consider the fundamental map
$$ F: X \times T \longrightarrow \bold P(\Cal Q) $$
and its associated fundamental involution
$$j: X \times T \longrightarrow X \times T$$
which have been studied in detail in the previous sections. We want to
apply theorem 6.13 to $F$ and $j$, therefore we must check wether the
assumptions (1),(2),(3),(4) of
this theorem are satisfied: \par \noindent
Assumption (3) is satisfied because $X$ is normal and $Pic(X)-Sing X \cong Pic
X $ (see [B1]3.1). By remark 2.18 assumption (4) is satisfied
too. Assumptions (1) and (2) concern the locus
$$
N(\Cal L) = \lbrace o \in X / \theta \text {is not an embedding at $o$} \rbrace
$$
\noindent a point $o \in N(\Cal L)$ satisfies them iff: \par \noindent
(1) $ o \times T \cap I_j$
has codimension $\geq 2$ in $T$, where $I_j$ is the indeterminacy locus of $j$
\par \noindentÊ (2)
let $p_1: X\times T \rightarrow X$ be the first projection, then $p_1 \cdot j$
is generically defined
and constant when restricted to $o \times T$. \par \noindent
\proclaim {(7.1) LEMMA}Ê \par \noindent
(i) Let  $o \in N(\Cal L) \cap X-Sing X$ then $o$ satisfies assumptions (1) and
(2) \par \noindent
(ii) let $o \in Sing X$ then $\theta^{-1}(o) = \lbrace o \rbrace$. \endproclaim
\demo {Proof} (i) At
first we show that condition (2) is satisfied by $o$. We distinguish two cases:
\par \noindent  CASE
(A) \it $\theta$ is not injective at $o$ \rm.  Then there exists $\overline o
\neq o$ such that
$\theta(o) = \theta(\overline o)$, we will assume $o = [\xi]$, $\overline {o} =
\overline {\xi}$.
Consider in $T$ the special closed set $I(o,\overline o)$ defined in 5.5 and
the complement $$
U_{o,\overline o} = T-I(o,\overline o) $$ \noindent of the special closed set
$I(o,\overline o)$
$=$ $I(o) \cup I(\overline o)$ defined in 5.5: we already know that,  putting
$$\Cal E = \xi(t),
V = H^0(\xi(t)) \quad \overline {\Cal E} = \overline {\xi(t)},\overline V =
H^0(\overlineÊ{\xi}(t))$$
and choosing $t$ in $U_{o,\overline o}$,  the dimension of$ V$, $\overline V$
is $4$
and moreover the pairs $(\Cal E,V)$, $(\overline {\Cal E}, \overline V)$ define
two quadrics $$
Q_t = Q(\Cal E,V), \quad {\overline Q}_t = Q(\overline {\Cal E},\overline V)$$
having rank $\geq 4$.
Actually, since $o$ is not in $Sing X$ and $t$ is not in $I(o)$, there is no
subline bundle $L$ of
$\Cal E$ with $h^0(L) \geq 2$ so that   $$ rank Q_t \geq 5 $$ \noindent by
prop. 1.11. Choosing
$t$  general in $U_{o,\overline o}$ (that is:  $t \in U_0-(U_0 \cap V_1(\xi)$)
$\xi(t)$ is
globally generated, hence  $$ Sing Q_t \cap C_t = \emptyset $$ \noindent by
prop. 1.9.  From
the definition of the fundamental map $F$ we have  $$ Q_t = F(o,t) \quad \text
{and}Ê\quad \overline
{Q}_t = F(\overline o,t)  $$ \noindent Since $F = \lambda \cdot (\theta \times
id_T)$ and $\theta(o) =
\theta(\overline o)$ it follows $F(o,t)$ = $F(\overline o,t)$ that is $$ Q_t =
\overline {Q}_t  $$
\noindent Since $Q_t$ has rank $\geq 5$ and $Sing Q_t \cap C_t = \emptyset$
there exist at most two
isomorphism classes of pairs defining $Q_t$, moreover they are either
isomorphic or dual (prop. 1.19,
lemma 1.18). In particular this holds for $(\Cal E,V)$, $(\overline {\Cal
E},\overline V)$: since $o
\neq \overline o$ $\xi$ is not isomorphic to $\overline {\xi}$ hence these two
pairs cannot be
isomorphic and they are dual. Then the rank of $Q_t$ is 6 and  $$ j(o,t) =
(\overline o,t)$$ \noindent
Since this holds for a general $t \in T$ $j$ is generically defined on $o
\times T$ and
$p_1 \cdot j:  o \times T \rightarrow X$ is the constant map onto the point
$\overline o$. Hence
condition (2) of theorem 6.13 is satisfied by $o$. The proof of the next case
is similar:  \par
\noindent CASE (B) \it the tangent map $d\theta_{o}: T_{X,o} \rightarrow
T_{\theta(X),\theta(o)}$ is
not an isomorphism. \rm Since we have already shown case (A) we can assume
$\theta^{-1}(o) = \lbrace
o \rbrace $: this clearly implies $d\theta_{o}$ not  injective. We consider as
above the quadric $Q_t
= F(o,t)$:  we know that its rank is 5 or 6 for general $t$ (precisely for $t
\in T-W(\xi)$). Let us
show that in this case the rank is 5:  since $F = \overline {\lambda} \cdot
(\theta \times id_T)$ and
$d\theta$ is not injective at $o$  $dF$ is not injective at $(o,t)$. We have
shown in proposition
3.18 that this is impossible if the rank of $Q_t$ is 6. Since the rank is 5
$(\Cal E,V)$ is
isomorphic to its dual so that  $j(o,t) = (o,t)$ for a general $t$. Hence $o$
satisfies conditions
(2) of theorem 6.13. To complete the proof of (i) we have now to show that $o$
satisfies also
condition (1). This is an immediate consequence of theorem 5.8 which says that,
with the prescribed
exceptions, $codim o \times T \cap I_j \geq 2$ if $p_1 \cdot j: o \times T
\rightarrow X$ is
constant. Indeed, since we are assuming $C$ not hyperelliptic, $g \geq 4$, $o$
non singular, only the
following exception is possible: \par \noindent there exists a double covering
$\pi: C \rightarrow Y$ of an elliptic curve and $o$ is the moduli point of the
bundle constructed
from $\pi$ as in the statement of theorem 5.8. Our claim is that $\theta$ is an
embedding at such a
special point $o$ so that  $o$ is not in $N(\Cal L)$: we will show this claim
at the end of the
paper.  \par \noindent (ii) Let $u \in X-Sing X$ we have already remarked that
$F(u,t)$ is a quadric
of rank $\geq 5$ if $t$ is general. On the other hand the rank of $F(o,t)$ is
always $\leq 4$ if $o =
[\xi]\in Sing X$: this is due to the fact that $\xi(t)$ has a subline bundle
$L$ with $h^0(L) \geq 2$
for every $t$, ($deg L = g+1$). Therefore $F(u,t) \neq F(o,t)$ at least for one
$t$ and hence $\theta
(u) \neq \theta (o)$. In particular $\theta^{-1}(o) \subset Sing X$, but it is
well known that
$\theta/Sing X$ is injective hence $\theta^{-1}(o) = o$. \enddemo \noindent
Finally we can complete
the \proclaim {PROOF OF THE MAIN THEOREM} \endproclaim  Assume there exists a
point $o \in X$ such
that \par \noindent (1) $\theta$ is not injective at $o$ in $X$ or \par
\noindent (2) $o$ is non
singular and $d\theta_o$ is not an isomorphism.  \par \noindent In case (1), by
statement (ii) of the
previous lemma 7.1, $o \in X-Sing X \cap N(\Cal L)$ hence, by statement (i) of
the same lemma, $o$
satisfies assumptions (1), (2) of theorem 6.13. In case (2) $o\in X-Sing X \cap
N(\Cal L)$ again
hence the same holds by the same lemma 7.1. Therefore we can apply to $o$ the
rigidity argument used
in theorem 6.13: since $o$ exists and satisfies assumptions (1) and (2) $j$ is
induced by a non
identical birational involution $f$ on $X$. That is  $$j = f \times id_T$$.
\noindent Once we have
$f$ we conclude by showing that $$\theta \cdot f = \theta$$ so that
$deg(\theta) \geq 2$. For this it
suffices to produce a subset $S$ of $X$ such that $\theta(S)$ spans the ambient
space $\mid 2\Theta
\mid$ of $\theta(X)$ and $f/S = id_S$ (see remark 6.14). Let $S$ = $Sing X$
then $\theta(S)$ is the
Kummer variety of $J$ naturally embedded in $\mid 2\Theta \mid$ and of course
it spans this space.
Let $s \in S$: since, for a general $t$, $j(s,t) = (s,t)$ it follows $f(s)=s$
and $f/S
= id_S$. Therefore $deg(\theta) \geq 2$: since $C$ is not hyperelliptic this is
impossible ([B1]),
hence the point $o$ cannot exists. This implies the theorem.  \enddemo \bigskip
\noindent COMPLETION
OF THE PROOF OF LEMMA 7.1 (i) (the bielliptic case)  Let $\pi: C \rightarrow Y$
be a double cover of
an elliptic curve, $\eta$ an irreducible rank two vector bundle on $Y$ of
degree 1: as claimed in
the proof of lemma 7.1 (i) we must show that $\theta$ is an embedding at the
moduli point $o$  of
$$\xi   \pi^*\eta \otimes L$$ \noindent where $L^2 \cong \omega_C \otimes det
\pi^*\eta^*$. We recall
that $\xi$ is stable (cfr.[LN]) and that, fixing its determinant, $\eta$ is
unique. Hence for every
$\pi$ there are $2^{2g}$ exceptional points like $o$: one for each $L$. We will
just sketch
proof: by [LN] $\xi$ has a family of subline bundles of maximal degree $g-2$
which is  exactly the
following $$ Y_L = \lbrace L(f), f \in \pi^*Pic^0(Y) \rbrace $$ \noindent Let
$C_t \subset \bold
P^{g+2}_t$ be the curve $C$ embedded by $\omega_C(2t)$, as usual we can choose
$t \in T$ so that in
the quadric $$ Q_t = F(o,t) $$
\noindent has rank $\geq 5$ and $C_t \cap Q_t = \emptyset$. We can also assume
$h^0(L(f+t)) = 1$
and $h^1(L(f-t)) = 3$ for each $L(f) \in Y_L$. Then, denoting by $< D_f >$ the
linear span in $\bold
P^{g+2}_t$ of the unique element $ D_f \in \mid L(f+t) \mid$, it follows
$$codim <D_f> = 3$$
\noindent and moreover
$$
<D_f> \subset Q_t
$$
\noindent because $h^0(\xi(t-D_f))Ê\geq 1$. Obviously $<D_f> \neq <D_g>$ if $f
\neq g$; moreover,
since they are maximal subspaces of $Q_t$ in the same ruling we have
$$
dim <D_f> \cap <D_g> = g-5
$$
\noindent Now observe that, if the rank of $Q_t$ is $5$, then $dim
Sing Q_t = g-5$ so that $<D_f> \cap <D_g>$ $=$ $Sing Q_t$. Assuming this and
considering the variety
$R = \cup <D_f>, \quad f \in \pi^*Pic^0(Y)$ it is not difficult to deduce that
$R$ would be a cone of
vertex $Sing Q_t$ over $\overline R \subset \overline Q \subset \bold P^4$,
where $\overline R$ is a
smooth scroll over $Y$ and $\overline Q$ a smooth quadric. Since this is
impossible $Q_t$ must have
rank 6. Then all the assumptions of proposition 3.18 are satisfied by the pair
$(o,t)$ and hence
$(dF_t)_o = (d\lambda)_o \cdot (d\theta)_o $ is injective. Therefore
$d\theta_o$ is injective. It
remains to show that $\theta^{-1}(o) = \lbrace o \rbrace$: essentially this
follows from the remark
that $\Theta_{\xi}$ contains the irreducible components  $$ W-Y_L \quad , \quad
Y_L-W \quad \text
{(where $W = C(g-2)$})$$ \par \noindent \noindent Assume $\theta([\overline
\xi]) = \theta(o)$: by
lemma 7.1 (ii) $\overline {\xi}$ is stable, the main point is to deduce from
$W-Y_L \subset
\Theta_{\overline {\xi}}$ ($ = \Theta_{\xi}$) that $L(f)$ is a subline bundle
of $\overline {\xi}$ for
every $L(f) \in Y_L$: we omit this. Since, by the quoted result of
Lange-Narasimhan ([LN]thm.5.1),
the latter property characterizes $\xi$ it follows  $\overline {\xi} \cong
\xi$. \par 

\bigskip
\noindent
AUTHORS' ADDRESS: \par \noindent
S. Brivio, Universit\'a di Torino, Dipartimento di Matematica,
 via Carlo Alberto 10, 1O123 TORINO Italy

\noindent
A. Verra, Universit\'a di Roma "La Sapienza", Dipartimento di Matematica,
Piazzale A. Moro 2,
00185 ROMA Italy.

\end